\documentclass[11pt, twoside]{amsart}

\usepackage[utf8]{inputenc}
\usepackage[T1]{fontenc}
\usepackage{amssymb,amsmath,amstext}
\usepackage{hyperref}

\newtheorem{theor}{Theorem}[section]
\newtheorem{lem}[theor]{Lemma}
\newtheorem{defin}[theor]{Definition}

\newtheorem{notation}[theor]{Notation}

\newtheorem{cor}[theor]{Corollary}
\newtheorem{rem}[theor]{Remark}

\newtheorem{fact}[theor]{Fact}
\newtheorem{claim}[theor]{Claim}

\newtheorem{assump}[theor]{Assumption}
\newtheorem{observation}[theor]{Observation}

\numberwithin{equation}{section}

\newcommand{\acl}{\mathrm{acl}}
\newcommand{\dcl}{\mathrm{dcl}}
\newcommand{\tp}{\mathrm{tp}}

\newcommand{\cl}{\mathrm{cl}}

\newcommand{\es}{\emptyset}

\newcommand{\su}{\mathrm{SU}}

\newcommand{\crd}{\mathrm{crd}}

\newcommand{\meq}{^{\mathrm{eq}}}

\newcommand{\mcC}{\mathcal{C}}

\newcommand{\mcG}{\mathcal{G}}

\newcommand{\mcK}{\mathcal{K}}

\newcommand{\mcM}{\mathcal{M}}
\newcommand{\mcN}{\mathcal{N}}

\newcommand{\mcP}{\mathcal{P}}

\newcommand{\mcR}{\mathcal{R}}

\newcommand{\mcU}{\mathcal{U}}

\newcommand{\mbR}{\mathbf{R}}

\newcommand{\mbbN}{\mathbb{N}}
\newcommand{\mbbR}{\mathbb{R}}

\newcommand{\ind}{\raisebox{-2pt}[5pt][0pt]{$\smile$} \hspace*{-6.8pt}\raisebox{3pt}[5pt][0pt]{$|$} \; \: }
\newcommand{\nind}{\raisebox{-2pt}[5pt][0pt]{$\smile$} 
\hspace*{-6.8pt}\raisebox{3pt}[5pt][0pt]{$|$}\hspace*{-6.8pt}
\raisebox{3pt}[5pt][0pt]{$\diagup$} }

\newcommand{\rng}{\mathrm{rng}}

\title[Binary simple homogeneous structures]
{Binary simple homogeneous structures}

\author{Vera Koponen}

\address{Vera Koponen, Department of Mathematics, Uppsala University, Box 480,
75106 Uppsala, Sweden.}

\email{vera.koponen@math.uu.se}

\date{7 September 2016}

\begin{document}

\begin{abstract}
We describe all binary simple homogeneous structures $\mcM$ 
in terms of $\es$-definable equivalence relations on $M$, which ``coordinatize'' $\mcM$ and control dividing, and
extension properties that respect these equivalence relations.

\end{abstract}

\maketitle

\section{Introduction}\label{Introduction}

\noindent
We describe the fine structure of binary simple homogeneous structures to the extent that seems
feasible without further assumptions and with known concepts and methods from infinite model theory.
In this respect, this article completes the earlier work on this topic by Aranda Lop\'{e}z \cite{AL}, Ahlman \cite{AK} and the 
present author \cite{AK, Kop16a, Kop-one-based, Kop16b}.
Before discussing the results, we explain what ``homogeneity'' means here, and give some background.

We call a structure $\mcM$ {\em homogeneous} if it is countable, has a finite relational vocabulary (also called signature)
and every isomorphism between finite substructures of $\mcM$ can be extended to an automorphism of $\mcM$.\footnote{
The expressions {\em finitely homogeneous} and {\em ultrahomogeneous} are also used for the same notion.}
For a countable structure $\mcM$ with finite relational vocabulary, being homogeneous is equivalent to having
elimination of quantifiers \cite[Corollary~7.42]{Hod}; it is also equivalent to being a {\em Fra\"{i}ss\'{e} limit}
of an {\em amalgamation class} of finite structures \cite{Fra54, Hod}.
A structure with a relational vocabulary will be called {\em binary} if every relation symbol is unary or binary.
Certain kinds of homogeneous structures have been classified.
This holds for homogeneous partial orders, graphs, directed graphs, finite 3-hypergraphs,
and coloured multipartite graphs \cite{Che98, Gar, GK, Lach84, LachTripp, LT, LW, Schm, Shee}.
For a survey about homogeneous structures, including their connections to permutation groups, Ramsey theory,
topological dynamics and constraint satisfaction problems, see \cite{Mac11} by Macpherson.

A detailed theory, due to Lachlan, Cherlin, Harrington, Knight and Shelah
\cite{CL, KL, Lach84, Lach97, LS}, 
also exists for stable infinite homogeneous structures, for any finite relational language,
which describes them in terms of (finitely many) dimensions and $\es$-definable indiscernible sets (which may live in $\mcM\meq$);
see \cite{Lach97} for a survey. This theory also sheds light on finite homogeneous structures.
But we seem to be a very long way from a classification of (even binary) finite homogeneous structures.
This has consequences for (eventual) classifications of infinite homogeneous structures, for the following reason.
Suppose that $\mcN$ is a finite (binary) homogeneous structure. Let $\mcM$ be the disjoint union of $\omega$ copies
of $\mcN$ and add an equivalence relation such that each equivalence class is exactly the set of elements
in some copy of $\mcN$. Then $\mcM$ is a (binary) stable
homogeneous structure. Hence a classification of all (binary) stable homogeneous structures presupposes an
equally detailed classification of all (binary) finite homogeneous structures.
Thus we ignore the inner structure of such (``very local'') finite ``blocks'' as the copies of $\mcN$ in the example, 
and focus on the ``global fine structure'' of an infinite structure $\mcM$.

The notion of simplicity generalizes stability and implies that there is a quite useful notion of independence.
Moreover, there are interesting (binary) simple homogeneous structures which are unstable, such as the Rado graph
and (other) homogeneous metric spaces with a finite distance set. (More about this is Section~\ref{Metric spaces}).
From this point of view it is natural, and seems feasible, to study simple homogeneous structures.
From now on when saying that a structure is simple we assume that it is infinite, so ``simple and homogeneous'' implies
that it is countably infinite.
The theory of binary simple homogeneous structures has similarities to the theory of stable homogeneous structures, but
also differences. Every stable (infinite) homogeneous structure is $\omega$-stable, hence superstable, with finite SU-rank.\footnote{
This is explained a little bit in the introduction of \cite{Kop16a}.
The SU-rank is usually called U-rank in the context of stable structures.}
Analogously, every binary simple homogeneous structure is supersimple with finite SU-rank (which is bounded by the number
of 2-types over $\es$) \cite{Kop16a}. However, the rank considered in the work on stable homogeneous structures is
Shelah's ``$CR( \ , 2)$-rank'' \cite[p. 55]{She}. 
This rank is finite for stable homogeneous structures, but it is infinite for the Rado graph.
If $\mcM$ is stable (infinite) and homogeneous and $C \subseteq M\meq$ is $\es$-definable and such that,
on $C$, there is no $\es$-definable nontrivial equivalence relation, then $C$ is an indiscernible set.
This is not true in general for (binary) simple homogeneous structures, as witnessed again by the Rado graph.

Suppose that $\mcM$ is binary, simple, and homogeneous.
We already mentioned that $Th(\mcM)$, the complete theory of $\mcM$, is supersimple with finite SU-rank.
It is also known that $Th(\mcM)$ is 1-based and has 
trivial dependence/forking~\cite[Fact~2.6 and Remark~6.6]{Kop16b}.
If $\mcM$ is, in addition, {\em primitive}, then $\mcM$ has SU-rank 1 and is a random structure~\cite{Kop16b}.
(See Section~\ref{Simple homogeneous structures}  for a definition of `primitive structure'.)
Before stating the main results of this article, we note that, although the definition (above) of `homogeneous structure' involves the
assumption that the structure is countable, the main results hold for {\em every} model of $Th(\mcM)$.
The reason is that, $\mcM$ (being homogeneous) is $\omega$-categorical and hence $\omega$-saturated.
So if elements could be found in some $\mcN \models Th(\mcM)$ such that one of the statements (a)--(d) below fails in $\mcN$, 
then such elements could also be found in $\mcM$.

\medskip

\noindent
{\bf Main results} (Theorems~\ref{theorem coordinatizing M by equivalence classes} 
and~\ref{main result regarding near randomness}).
{\em
Suppose that $\mcM$ is binary, simple, and homogeneous (hence supersimple with finite SU-rank and trivial dependence).
Let $\mbR$ be the (finite) set of all $\es$-definable equivalence relations on $M$.
If $a \in M$ and $R \in \mbR$, then $a_R$ denotes the $R$-equivalence class of $a$ as an element of $M\meq$.
\begin{itemize}
\item[(a)] {\rm Coordinatization by equivalence relations:} 
For every $a \in M$, if $\su(a) = k$, then there are $R_1, \ldots, R_k \in \mbR$,
depending only on $\tp(a)$,  such that
$a \in \acl(a_{R_k})$, $\su(a_{R_1}) =~1$, 
$R_{i+1} \subset R_i$ and 
$\su(a_{R_{i+1}} / a_{R_i}) = 1$ for all $1 \leq i < k$
(or equivalently, $\su(a / a_{R_i}) = k-i$ for all $1 \leq i \leq k$).

\item[(b)] {\rm Characterization of dividing:} 
Suppose that $a, b, \bar{c} \in M$ and $a \underset{\bar{c}}{\nind} b$.
Then there is $R \in \mbR$ such that $a \underset{\bar{c}}{\nind} a_R$ and $a_R \in \acl(b)$
(and thus $a_R \notin \acl(\bar{c})$).

\item[(c)] {\rm Characterization of dividing in the symmetric case:} 
Suppose that all binary $\es$-definable relations on $M$ are symmetric.
If $a, b, \bar{c} \in M$ and $a \underset{\bar{c}}{\nind} b$,
then  there is $R \in \mbR$ such that $a \underset{\bar{c}}{\nind} a_R$ and $R(a, b)$ (hence $a_R \in \acl(b)$,
$a_R \notin \acl(\bar{c})$ and thus $\neg E(a, c)$ for every $c \in \bar{c}$).

\item[(d)] {\rm Extension properties:}
Let $a, b, c, \bar{d} \in M$.\\
(i) There is $R \in \mbR$ such that
$c \underset{c_R}{\ind} \bar{d}$.\\
(ii) If for some $R$ as in part~(i),
\begin{itemize}
\item[] $a \underset{c_R}{\ind} c$, \ $b \underset{c_R}{\ind} \bar{d}$ \ and \
$\tp(a / \acl(c_R)) = \tp(b / \acl(c_R))$, where `$\acl$' is taken in~$M\meq$,
\end{itemize}
then there is $e \in M$ such that $\tp(e, c) = \tp(a, c)$ and $\tp(e, \bar{d}) = \tp(b, \bar{d})$.
Otherwise such $e$ may not exist (in any elementary extension of $\mcM$), not even when $\bar{d}$ is a single element.
\end{itemize}
}

\noindent
In parts~(b) and~(c) we only consider singletons $a$ and $b$ because $\mcM$ has trivial dependence.
We will show (in Section~\ref{Bipedes with bicoloured legs}) that the ``symmetry condition'' in part~(c) cannot be removed;
in other words, the conclusion in part~(b) cannot be strengthened so that it becomes identical to the conclusion in part~(c).
Further remarks on (a)--(c) are made in Remark~\ref{remarks on the coordinatization by equivalence relations}.
Regarding part~(d)(ii),
the conditions that $a \underset{c_R}{\ind} c$, $b \underset{c_R}{\ind} \bar{d}$ and 
$\tp(a / \acl(c_R)) = \tp(b / \acl(c_R))$, are just the premisses (in the present context) 
of the independence theorem for simple theories. 
So the interesting part, with respect to~(d)(ii), is that if (for every $R$ as in~(i)) these premisses are not satisfied,
then a ``common extension'' may not exist. 
Thus we do not, in general, get anything ``for free'' beyond what the independence theorem guarantees.
From this, one may get the impression that common extensions of types like in~(d) are unusual.
But note that, by part~(i) of~(d), we can always find a $\es$-definable equivalence relation $R$ such that
$c \underset{c_R}{\ind} \bar{d}$.
Therefore I would say that (by part~(ii) of~(d)), in a binary simple structure, common extensions of two types do exist
as long as we respect all $\es$-definable equivalence relations and some other ``reasonable'' 
conditions related to them.
The examples in sections~\ref{Cross cutting equivalence relations} --~\ref{omega-pedes} show that these conditions
are, in fact, necessary.
The reason that~(d) only considers an extension of two 1-types (one of which has only one parameter $c$) is that,
since $\mcM$ is binary with elimination of quantifiers,
the problem of extending more than two $k$-types (with finite parameter sets) can be reduced to a 
finite sequence of ``extension problems'', 
each of which involves only two 1-types and one of the types has only one parameter.
More about this is said in the beginning of Section~\ref{Extension properties}.

From the proofs of the main results, one can extract information about 
$\omega$-categorical (not necessarily binary or homogeneous) supersimple structures with finite SU-rank and trivial dependence.
This information is presented in Corollaries~\ref{information about nonbinary structures}
and~\ref{information about nonbinary primitive structures}, and may be useful in future studies of
nonbinary simple homogeneous structures.

Now we turn to problems about simple homogeneous structures.
If $\mcM$ is stable and homogeneous, then $\mcM$ has the {\em finite submodel property}\footnote{
This means that every sentence which is satisfied by $\mcM$ is also satisfied by a finite substructure of $\mcM$.}
and $Th(\mcM)$ is decidable. (For the first result, see \cite[Proposition~5.1]{Lach97} or
\cite[Lemma~7.1]{KL}; for the second, see the proof of Theorem~5.2 in \cite{Lach97}.)
It is still not settled whether every binary simple homogeneous structure has the finite submodel property, nor whether its
theory must be decidable. But my guess is that the answer is `yes' to both questions.\footnote{
I believe that a positive answer may involve some probabilistic arguments in the style of the 0-1 law for finite graphs/structures
\cite[Lemma~7.4.6]{Hod}, 
or (more generally) the arguments by Ahlman in \cite{Ahl}.
But the appearance of multiple $\es$-definable equivalence relations, which may refine or ``cut'' each other 
in various ways, seems to make a proof of this kind a bit more than a ``straightforward generalization'' of 
known arguments.}

Regarding nonbinary simple homogeneous structures, I would say that all core problems are unsolved.
The answer is unknown to each of these questions, where we assume that $\mcM$ is (nonbinary) simple and homogeneous:
Must $Th(\mcM)$ be supersimple? If $Th(\mcM)$ is supersimple, must it have finite SU-rank?.
Must $Th(\mcM)$ be 1-based? Must $Th(\mcM)$ have trivial dependence?
(If $\mcM$ is supersimple, the last two problems are tightly connecteds to the problem of which kinds of definable pregeometries,
induced by algebraic closure, there can be on the realizations, in $M\meq$, of types of SU-rank~1.)
If $\mcM$ is supersimple with SU-rank 1, what possibilities are there for the fine structure of $\mcM$
(according to some ``reasonably'' informative classification)?
Even if we add `primitivity' and `trivial dependence' to the assumptions of the last question, the answer is unknown.

Here follows an outline of the article.
Section~\ref{Preliminaries} explains the notation and terminology that will be used, and
gives background regarding 
homogeneous (or just $\omega$-categorical) simple structures.
Section~\ref{Coordinatization} describes the ``coordinatization'' developed in \cite[Section~3]{Djo06}
for $\omega$-categorical, supersimple structures with finite SU-rank and trivial dependence
(or equivalently, $\omega$-categorical simple 1-based structures with trivial dependence).
This coordinatization will be the framework in Sections~\ref{The main technical lemmas} 
and~\ref{Coordinatization by equivalence relations}.
In Section~\ref{The main technical lemmas} we prove the main technical lemmas, on which the main results rest.
In Section~\ref{Coordinatization by equivalence relations} we prove (a)--(c) from the main results above.
(This involves proving that every ``coordinate'' in the sense of Section~\ref{Coordinatization} is interalgebraic with
a new coordinate $a_R$ where $a \in M$ and $R$ is a $\es$-definable equivalence relation on $M$.)
In Section~\ref{Extension properties} we partially prove part~(d) above, with the help of part~(b).
To complete the proof of~(d), we also need to construct ``counterexamples'', which is done
in Sections~\ref{Cross cutting equivalence relations} --~\ref{omega-pedes}.
Section~\ref{Metric spaces} is an exposition of results by Conant \cite{Con-neostab}
about homogeneous metric spaces, which concretize the main results of this article in that context.

\section{Preliminaries}\label{Preliminaries}

\subsection{Notation and terminology}

Structures will be denoted by calligraphic letters, usually $\mcM$ or $\mcN$ in which case their universes are denoted $M$ or $N$,
respectively. Finite sequences (and {\em only finite} sequences) are denoted by $\bar{a}, \bar{b}, \ldots, \bar{x}, \bar{y}, \ldots$.
The concatenation of $\bar{a}$ and $\bar{b}$ is denoted $\bar{a}\bar{b}$, but sometimes we also write $(\bar{a}, \bar{b})$
(like when using the type notation $\tp(\bar{a}, \bar{b})$). The set of elements that occur in $\bar{a}$
(in other words, the range/image of $\bar{a}$) is denoted $\rng(\bar{a})$.
But when the order of $\bar{a}$ does not matter,
we often abuse notation and (notationally) identify the sequence $\bar{a}$ with the set $\rng(\bar{a})$.
So we may write things like `$a \in \bar{a}$' instead of '$a \in \rng(\bar{a})$'.
When $a$, $b$ and $c$ are single elements we sometimes write `$ab$' for the pair `$(a, b)$', or
`$abc$' for the triple `$(a, b, c)$', and similarly for longer tuples.
Further, we often write `$\bar{a} \in A$' when meaning that $\bar{a}$ is a finite sequence such that $\rng(\bar{a}) \in A$.
If we may emphasize that the length of $\bar{a}$ (denoted $|\bar{a}|$) is $n$, then we may write $\bar{a} \in A^n$.

As usual, `$\acl_\mcM$', `$\dcl_\mcM$', and  `$\tp_\mcM$' denote the algebraic closure, definable closure,
and type (of a set or sequence) in the structure $\mcM$; and if $A \subseteq M$, then 
$S_n^\mcM(A)$ is the set of $n$-types over $A$ with
respect to $Th(\mcM)$, the complete theory of $\mcM$.
The notation `$\bar{a} \equiv_\mcM \bar{b}$' means the same as
`$\tp_\mcM(\bar{a}) = \tp_\mcM(\bar{b})$'.
The notation `$\bar{a} \equiv^{at}_\mcM \bar{b}$' means that $\bar{a}$ and $\bar{b}$
satisfy exactly the same atomic formulas with respect to $\mcM$.
In sections~\ref{Coordinatization} --~\ref{Extension properties} the structure $\mcM$ is fixed and we work in $\mcM\meq$,
so for brevity we will, in those sections, omit the subscript `$\mcM\meq$' and write for example `$\tp$' instead of  `$\tp_{\mcM\meq}$'.
We remind about this again in Notation~\ref{remark about acl and dcl}.
If $p(\bar{x})$ is a type (or formula), then $p(\mcM)$ denotes the set of realizations of $p$ in $\mcM$.

If $R$ is a $\es$-definable equivalence relation on $M^n$ for some $n < \omega$, then we may also call $R$ a {\em sort}.
For every such $R$ and every  $\bar{a} \in M^n$, $[\bar{a}]_R$ denotes the $R$-equivalence class of~$\bar{a}$.
When we view $[\bar{a}]_R$ as an {\em element} of $M\meq$ we write $\bar{a}_R$ to emphasize this.
If $A \subseteq M\meq$ then we say that {\em only finitely many types are represented in $A$}
if there are only finitely many sorts $R$ such that for some
$n < \omega$ and $\bar{a} \in M^n$, $\bar{a}_R \in A$.

When saying that $\mcM$ is $\omega$-categorical, (super)simple, 1-based, or that $\mcM$ has finite SU-rank, then we mean that
$Th(\mcM)$ is $\omega$-categorical,
(super)simple, 1-based, or that $Th(\mcM)$ has finite SU-rank, respectively.

A {\em pregeometry} (or {\em matroid}) is a pair $(X, \cl)$ where $X$ is a set and 
$\cl : \mcP(X) \to \mcP(X)$ satisfies certain conditions (see \cite[Chapter~4.6]{Hod}).
We say that a pregeometry $(X, \cl)$ is {\em trivial} if for all $Y  \subseteq X$,
$\cl(Y) = \bigcup_{a \in Y}\cl(\{a\})$.

\subsection{$\omega$-Categorical structures}\label{omega-categorical and simple structures}

Since homogeneous structures have elimination of quantifiers, it follows from the 
well-known characterization of $\omega$-categoricity \cite[Theorem~7.3.1]{Hod}, that every 
infinite homogeneous structure is $\omega$-categorical.
We now state some basic facts about $\mcM\meq$ when $\mcM$ is $\omega$-categorical.
These will tacitly be used throughout the article.

\begin{fact}\label{fact about isolated types}
Suppose that $\mcM$ is $\omega$-categorical and assume that only finitely many sorts are
represented in $A \subseteq M\meq$.
\begin{itemize}
\item[(i)] For every $n < \omega$ and every finite $B \subseteq M\meq$, only finitely many types from
$S_n^{\mcM\meq}(\acl_{\mcM\meq}(B))$ are realized by tuples in $A^n$.
\item[(ii)] For every finite $B \subseteq M\meq$, $A \cap \acl_{\mcM\meq}(B)$ is finite.
\item[(iii)] For every $\bar{a} \in M\meq$ and every finite $B \subseteq M\meq$, the types
$\tp_{\mcM\meq}(\bar{a} / B)$ and 
$\tp_{\mcM\meq}(\bar{a} / \acl_{\mcM\meq}(B))$ are isolated.
\end{itemize}
\end{fact}

\noindent
For some explanations of the above claims, see \cite[Section~2.4]{AK}.
Part~(iii) of Fact~\ref{fact about isolated types} will usually not be used in the form stated above, but rather we use
the following (namely {\em $\omega$-homogeneity} and a variant of it), which are proved  straightforwardly from 
Fact~\ref{fact about isolated types}~(iii):

\begin{fact}\label{generalized omega-homogeneity} 
Suppose that $\mcM$ is $\omega$-categorical.
\begin{itemize}
\item[(i)]  If $\bar{a}, \bar{b}, c \in M\meq$ and $\bar{a} \equiv_{\mcM\meq} \bar{b}$, then there is $d \in M\meq$ such that 
$\bar{a}c \equiv_{\mcM\meq} \bar{b}d$.
\item[(ii)] If $\bar{a}, \bar{b}, \bar{c}, \bar{e} \in M\meq$ and 
\[
\tp_{\mcM\meq}\big(\bar{a} / \acl_{\mcM\meq}(\bar{e})\big) \ = \ 
\tp_{\mcM\meq}\big(\bar{b} / \acl_{\mcM\meq}(\bar{e})\big),
\]
then there is $\bar{d} \in M\meq$ such that 
\[
\tp_{\mcM\meq}\big(\bar{a}\bar{c} / \acl_{\mcM\meq}(\bar{e})\big) \ = \ 
\tp_{\mcM\meq}\big(\bar{b}\bar{d} / \acl_{\mcM\meq}(\bar{e})\big).
\]
\end{itemize}
\end{fact}

\subsection{Simple homogeneous structures}\label{Simple homogeneous structures}

\noindent
We assume basic knowledge about simple structures as can be found in~\cite{Wag}, for instance,
but nevertheless recall a couple of things.
When saying that a structure is simple we automatically assume that it is infinite.\footnote{
Thus we do not follow the terminology of the work on stable homogeneous structures, where every finite structure is 
considered to be stable.}
Since $\omega$-categorical simple theories have elimination of hyperimaginaries \cite[Corollary~6.1.11]{Wag},
the independence theorem of simple theories 
\cite[Theorem~2.5.20]{Wag} takes the following form if the involved sets of parameters are finite and 
$\mcM$ is $\omega$-categorical and simple:
\begin{itemize}
\item[] Suppose that $\bar{a}, \bar{b} \in M\meq$, $A, B, C \subseteq M\meq$ are finite, 
$\bar{a} \underset{C}{\ind} A$, $\bar{b} \underset{C}{\ind} B$,
and
\[
\tp_{\mcM\meq}\big(\bar{a} / \acl_{\mcM\meq}(C)\big) \ = \
\tp_{\mcM\meq}\big(\bar{b} / \acl_{\mcM\meq}(C)\big).
\]
Then there is $\bar{d} \in M\meq$ such that 
\begin{align*}
\tp_{\mcM\meq}\big(\bar{d} / A \cup \acl_{\mcM\meq}(C)\big) \ &= \
\tp_{\mcM\meq}\big(\bar{a} / A \cup \acl_{\mcM\meq}(C)\big) \ \text{ and } \\
\tp_{\mcM\meq}\big(\bar{d} / B \cup \acl_{\mcM\meq}(C)\big) \ &= \
\tp_{\mcM\meq}\big(\bar{b} / B \cup \acl_{\mcM\meq}(C)\big).
\end{align*}
\end{itemize}

Note that if $\mcM$ is $\omega$-categorical and supersimple with finite SU-rank, then 
(since $S_1^\mcM(\es)$ is finite) there is $n < \omega$ such that
$\su(p) \leq n$ for every $S_1^\mcM(\es)$.
Before recalling what is known from before about binary simple homogeneous structures, we give the 
definition of trivial dependence (also called `trivial forking', or `totally trivial' in \cite{Goode}).

\begin{defin}{\rm
A simple complete theory $T$ has {\em trivial dependence} if 
for all $\mcM \models T$ and all $A, B, C \subseteq M\meq$, if $A \underset{C}{\nind} B$, then
$A \underset{C}{\nind} b$ for some $b \in B$.
We say that a simple structure $\mcM$ has {\em trivial dependence} if $Th(\mcM)$ has it.
}\end{defin}

\begin{fact}\label{facts about binary simple homogeneous structures}
Suppose that $\mcM$ is binary, simple, and homogeneous. Then:
\begin{itemize}
\item[(i)] $\mcM$ is supersimple with finite SU-rank (which is bounded by the number of complete 2-types over $\es$).
\item[(ii)] $\mcM$ has trivial dependence.
\item[(iii)] $\mcM$ is 1-based.
\end{itemize}
\end{fact}

\noindent
Part~(i) is given by \cite[Theorem~1]{Kop16a}.
Parts~(ii) and~(iii) are consequences of \cite[Corollary~6]{Kop16a}, \cite[Lemma~1]{Goode}, \cite[Corollary~4.7]{HKP},
\cite[Corollary~3.23]{DK} and \cite[Theorem~1.1]{Mac91}; this is explained in more detail
in the text surrounding Fact~2.6 in \cite{Kop16b} and in Remark~6.6 of the same article.

We call a structure $\mcM$ {\em primitive} if there there is no nontrivial $\es$-definable equivalence relation on $M$
(where by nontrivial we mean that there are at least two equivalence classes and at least one equivalence class has at
least two elements).

\begin{fact}\label{facts about binary primitive simple homogeneous structures}
Suppose that $\mcM$ is binary, primitive, simple and homogeneous. Then:
\begin{itemize}
\item[(i)] $\mcM$ has SU-rank 1.
\item[(ii)] $\mcM$ is a random structure in the sense of \cite[Definition~2.1]{Kop16b}.
\end{itemize}
\end{fact}

\noindent
Part~(i) is given by \cite[Theorem~1.1]{Kop16b}. 
Part~(ii) is a consequence of part~(i) and \cite[Proposition~3.3.3]{AL}, where the later result says that every
binary simple homogeneous structure of SU-rank 1 is a random structure.
From Theorem~\ref{theorem coordinatizing M by equivalence classes}~(i)
(i.e. part~(a) of the `main results' in the introduction), it follows that
part~(i) of Fact~\ref{facts about binary primitive simple homogeneous structures}
still holds if the assumption about `primitivity' is replaced with the condition that there is no 
$\es$-definable equivalence relation on $M$ which has infinitely many infinite equivalence classes.

Fact~\ref{facts about binary primitive simple homogeneous structures}~(i) fails without the binarity condition 
as shown by Example~2.7 in \cite{Kop16b}, which is primitive, homogeneous, and superstable with SU-rank 2 (but nonbinary).
It is also {\em not} a random structure.
Consequently also part~(ii) of Fact~\ref{facts about binary primitive simple homogeneous structures}
fails without the binarity condition.
But in fact it fails (without the binarity condition) in a stronger sense.
Because the generic tetrahedron-free 3-hypergraph is primitive, homogeneous, supersimple with SU-rank 1 and 1-based,
but not a random structure. All mentioned properties of the generic tetrahedron-free 3-hypergraph, except for the 1-basedness, 
have been known for a long time. Results which imply that it is 1-based were recently proved by  Conant~\cite{Con}
and by the present author~\cite{Kop16b}.

\section{Coordinatization}\label{Coordinatization}

\noindent
Throughout this section we assume that $\mcM$ is $\omega$-categorical, supersimple with finite SU-rank
and trivial dependence (hence it is 1-based).
Then the ``coordinatization'' results of Section~3 in \cite{Djo06} apply to $\mcM$.
We will now go through these results, since they are the framework in which the arguments
of sections~\ref{The main technical lemmas}
--~\ref{Extension properties}
take place.

\begin{notation}\label{remark about acl and dcl}{\rm
In this section and Sections~\ref{The main technical lemmas} --~\ref{Extension properties}, 
`$\tp$', `$\equiv$', `$\acl$', and `$\dcl$' will abbreviate 
`$\tp_{\mcM\meq}$', `$\equiv_{\mcM\meq}$', `$\acl_{\mcM\meq}$',
and `$\dcl_{\mcM\meq}$', respectively.
}\end{notation}

\begin{fact}\label{existence of coordinates}
Let $U \subseteq M\meq$ and suppose that only finitely many sorts are represented in $U$.
Then there are $0 < r < \omega$ and 
\[
C_0 \subseteq C_1 \subseteq \ldots \subseteq C_r \subseteq C \subseteq M\meq
\]
such that:
\begin{itemize}
\item[(i)] $U \subseteq C$, only finitely many sorts are represented in $C$, and $C$ is 
{\rm self-coordinatized} in the sense of \cite[Definition~3.3]{Djo06}.
\item[(ii)] $C$ and $C_i$ are $\es$-definable, for every $i = 1, \ldots, r$.
\item[(iii)] $C_0 = \es$ and, for every $n < h$ and every $c \in C_{n+1}$, $\su(c / C_n) = 1$.
\item[(iv)] $C \subseteq \acl(C_r)$.
\item[(v)] For every $1 < n \leq r$ and every $c \in C_n$, $\acl(c) \cap C_{n-1} \neq \es$.
\end{itemize}
\end{fact}

\begin{assump}\label{assumption about C}{\rm
In the rest of this section we suppose the following:
\begin{itemize}
\item[(a)] $M \subseteq U \subseteq M\meq$ and only finitely many sorts are represented in $U$.
\item[(b)] $C$ and $C_i$, for $i = 0, \ldots, r$, are as in Fact~\ref{existence of coordinates}.
\end{itemize}
}\end{assump}

\noindent
We can think of $C_r$ a set {\em coordinates of $C$} (and hence of $M$) and we call $h$ the
{\em height} of the coordinatization.

\begin{defin}\label{definition of crd}{\rm
(i) For every $\bar{c} \in C$ and every 
$0 \leq s \leq r$, let $\crd_s(\bar{c}) = \acl(\bar{c}) \cap C_s$. \\
(ii) We abbreviate `$\crd_r$' with `$\crd$'.
}\end{defin}

\noindent
Observe that for every $\bar{c} \in C$, $\crd(\bar{c})$ is finite.
We can think of $\crd(\bar{c})$ as the {\em coordinates} of $\bar{c}$ (with respect to the given coordinatization $C_r$)
and $\crd_s(\bar{c})$ as the {\em coordinates of $\bar{c}$ up to ``level'' $s$}.

\begin{fact}\label{fact about C and crd}
\begin{itemize}
\item[(i)] If $c \in C_r$, $d_1, \ldots, d_n \in M\meq$ and $c \in \acl(d_1, \ldots, d_n)$, 
then $c \in \acl(d_i)$ for some $1 \leq i \leq n$.\footnote{
This is \cite[Lemma~3.16]{Djo06}.}
\item[(ii)] For every $0 < s \leq r$, $(C_s \setminus C_{s-1}, \cl)$, where $\cl(A) = \acl(A) \cap (C_s \setminus C_{s-1})$ for
all $A \subseteq C_s \setminus C_{s-1}$, is a trivial pregeometry.\footnote{
This is an immediate consequence of \cite[Lemma~3.18]{Djo06}, because $C_s \setminus C_{s-1}$ is a $\es$-definable set and a
subset of the ($\es$-definable) set $N_s$ considered there \cite[Construction~3.13]{Djo06}.}
\item[(iii)] For every $\bar{c} \in C$
and every $0 \leq s \leq r$, $\crd_s(\bar{c}) = \bigcup_{c \in \rng(\bar{c})} \crd_s(c)$.
Thus the same holds for `$\crd$' in place of `$\crd_s$'.\footnote{
This is \cite[Lemma~5.4]{Kop09}.}
\item[(iv)] For all $\bar{c} \in C$, $\acl(\bar{c}) = \acl(\crd(\bar{c}))$.\footnote{
By definition, $\crd(\bar{c}) \subseteq \acl(\bar{c})$, so it suffices to show that
$\bar{c} \in \acl(\crd(\bar{c}))$.
By \cite[Lemma~5.1]{Djo06}, for every $c \in \bar{c}$, $c \in \acl(\crd(c))$. 
Thus the conclusion follows from part~(iii).}
\item[(v)] For all $\bar{a}, \bar{b} \in C$, $\bar{a}$ is independent from $\bar{b}$ over $\crd(\bar{a}) \cap \crd(\bar{b})$.\footnote{
Let $\bar{c}$ enumerate $\crd(\bar{a})$ and let $\bar{d}$ enumerate $\crd(\bar{b})$.
By part~(iv), $\acl(\bar{a}) = \acl(\bar{c})$ and $\acl(\bar{b}) = \acl(\bar{d})$, so
$\acl(\crd(\bar{a}) \cap \crd(\bar{b})) = \acl(\crd(\bar{c}) \cap \crd(\bar{d}))$.
Therefore it suffices to prove that $\bar{c}$ is independent from $\bar{d}$ over 
$\crd(\bar{c}) \cap \crd(\bar{d})$. Since $\bar{c}, \bar{d} \in C_r$, this is exactly the content of 
\cite[Lemma~5.16]{Kop09}.}
\end{itemize}
\end{fact}

\noindent
We note the following strengthening of part~(iii) of Fact~\ref{existence of coordinates}:

\begin{fact}\label{rank 1 over its coordinates on the previous level}
Let $0 \leq n < r$. For every $c \in C_{n+1} \setminus C_n$, $\su(c / \crd_n(c)) = 1$.
\end{fact}

\noindent
{\bf Proof.}
Suppose that $c \in C_{n+1} \setminus C_n$.
By Fact~\ref{existence of coordinates}~(iii), $\su(c / C_n) = 1$.
By supersimplicity, there is $\bar{d} \in C_n$ such that $\su(c / \bar{d}) = 1$.
Fact~\ref{fact about C and crd}~(v) implies that 
$c$ is independent from $\bar{d}$ over $\crd(c) \cap \crd(\bar{d})$,
so $\su(c / \crd(c) \cap \crd(\bar{d})) = 1$.
Since $\bar{d} \in C_n$ it follows from Fact~\ref{existence of coordinates}~(iii) that
$\crd(\bar{d}) \subseteq C_n$. Therefore $\su(c / \crd_n(c)) = 1$.
\hfill $\square$
\\

\noindent
The following generalization of Fact~\ref{fact about C and crd}~(v)
will be convenient to use.

\begin{lem}\label{a more general independence criterium}
Suppose that $\bar{a}, \bar{b}, \bar{c} \in C$. 
Then $\bar{a} \underset{\bar{c}}{\ind} \bar{b}$ if and only if 
$\crd(\bar{a}) \cap \crd(\bar{b}) \subseteq \acl(\bar{c})$.
\end{lem}

\noindent
{\bf Proof.}
Suppose that $\crd(\bar{a}) \cap \crd(\bar{b}) \subseteq \acl(\bar{c})$.
By extending the sequence $\bar{c}$ with new elements from $\crd(\bar{a}) \cap \crd(\bar{b})$, if necessary,
we may assume that $\crd(\bar{a}) \cap \crd(\bar{b}) \subseteq \bar{c}$.
By Fact~\ref{fact about C and crd}~(iii), $\crd(\bar{b}\bar{c}) = \crd(\bar{b}) \cup \crd(\bar{c})$,
so by Fact~\ref{fact about C and crd}~(v), 
$\bar{a}$ is independent from $\bar{b}\bar{c}$ over
\[
\crd(\bar{a}) \ \cap \ \big( \crd(\bar{b}) \cup \crd(\bar{c}) \big) \ = \
\big( \crd(\bar{a}) \cap \crd(\bar{b}) \big) \ \cup \  \big( \crd(\bar{a}) \cap \crd(\bar{c}) \big).
\]
So by monotonicity and the assumption that 
$\crd(\bar{a}) \cap \crd(\bar{b}) \subseteq \bar{c}$, it follows that
$\bar{a}$ is independent from $\bar{b}\bar{c}$ over $\bar{c}$.
Hence $\bar{a}$ is independent from $\bar{b}$ over $\bar{c}$.

Now suppose that $\bar{a} \underset{\bar{c}}{\ind} \bar{b}$.
For a contradiction, suppose that $d \in \crd(\bar{a}) \cap \crd(\bar{b})$ and $d \notin \acl(\bar{c})$.
Then $\su(d / \bar{c}) \geq 1$.
Using that $\acl(\bar{a}) = \acl(\crd(\bar{a}))$ (by Fact~\ref{fact about C and crd}~(iv),
we get, by the Lascar equation,
\[
\su(\bar{a} / \bar{c}) \ = \ \su(\bar{a}d / \bar{c}) \ = \ \su(\bar{a} / d\bar{c}) \ + \ \su(d / \bar{c}).
\]
Hence $\su(\bar{a} / d\bar{c}) < \su(\bar{a} / \bar{c})$. Therefore $\bar{a} \underset{\bar{c}}{\nind} d$,
and as $d \in \crd(\bar{b}) \subseteq \acl(\bar{b})$, we get $\bar{a} \underset{\bar{c}}{\nind} \bar{b}$.
\hfill $\square$

\begin{defin}\label{definition of E-s}{\rm
For every $0 \leq s \leq r$ and all $a, b \in C$, let
\[
E_s(a, b) \ \Longleftrightarrow \ \crd_s(a) = \crd_s(b) \ \text{ and } \
\tp\big( a \big/ \acl(\crd_s(a))\big) = \tp\big( b \big/ \acl(\crd_s(b))\big).
\]
}\end{defin}

\noindent
From 
Fact~\ref{fact about isolated types}~(iii) it is straightforward to derive the following:

\begin{fact}\label{fact about E-s}
For every $0 \leq s \leq r$, $E_s$ is a $\es$-definable equivalence relation on $C$.
\end{fact}

\begin{lem}\label{if c belongs to crd(a) and crd(b) then a and b have the same type over acl(c)}
We may, without loss of generality, assume that $C_r$ has the following property: 
for all $a, b \in M$ and all $c \in C_r$, if $ac \equiv bc$ then $\tp(a / \acl(c)) = \tp(b / \acl(c))$.
\end{lem}

\noindent
{\bf Proof.}
Let $c \in C_r$ and $q(x) = \tp(c)$.
Suppose that there are $a, b \in M$ such that  $\tp(a / \acl(c)) \neq \tp(b / \acl(c))$.
By Fact~\ref{fact about isolated types}~(i), only finitely many complete types over $\acl(c)$ are realized in $M$.
By part~(iii) of the same fact, each such type is isolated.
Let $p_1,  \ldots, p_n$ enumerate all complete 1-types over $\acl(c)$ which are realized in $M$.
For each $i$, choose a formula that isolates $p_i$ and let $\bar{d}_i$ be the parameters
(from $\acl(c)$)  that occur in that formula.
Let $\bar{d} = c\bar{d}_1 \ldots \bar{d}_n$.
Then $\acl(\bar{d}) = \acl(c)$.
As $\mcM\meq$ has elimination of imaginaries, there is $d \in M\meq$ such that $\dcl(d) = \dcl(\bar{d})$.
Let $q' = \tp(d)$.
Now remove from $C$ all $c' \in C$ which realize $q$
and then add to what is left of $C$ all $d' \in M\meq$ which realize $q'$.
Then the modified $C$ has the property that whenever $a, b \in M$, $c \in C$, $c$ realizes $q'$
and $ac \equiv bc$, then $\tp(a / \acl(c)) = \tp(b / \acl(c))$.
Since (by Assumption~\ref{assumption about C}) only finitely many types over $\es$ are realized in $C$, 
it follows that we can continue this procedure in finitely many steps and get (new) $C$ and $C_r \subseteq C$ such
that the conclusion of the lemma holds.
Since the types $q$ and $q'$ above are isolated and
every change of element in this process, say from $c$ to $d$, is such that $\acl(c) = \acl(d)$,
it follows that the new $C$ and $C_0 \subseteq \ldots \subseteq C_r$ that we get have all 
the properties of the earlier facts and lemmas in this section.
\hfill $\square$

\section{The main technical lemmas}\label{The main technical lemmas}

\noindent
Throughout this section we assume that $\mcM$ is binary, simple, and homogeneous.
By Fact~\ref{facts about binary simple homogeneous structures},
$\mcM$ is supersimple, 1-based, with finite SU-rank and with trivial dependence.
We thus adopt Assumption~\ref{assumption about C},
as well as Notation~\ref{remark about acl and dcl}.
However, the assumption that $\mcM$ is binary and homogeneous 
(as opposed to only $\omega$-categorical) is only used once at the end of the proof of
Lemma~\ref{different coordinates have different types over acl(emptyset)}
and once at the end of the proof of 
Lemma~\ref{when there are independent coordinates with the same type over a}.

The goal of this section is to prove the following:
\begin{itemize}
\item[] For all $0 < s \leq r$, $a \in M$ and $c_1, c_2 \in \crd_s(a) \setminus C_{s-1}$, if $c_1 \underset{\crd_{s-1}(a)}{\ind} c_2$ 
then $ac_1 \not\equiv ac_2$.
\end{itemize}
This is also the statement of Lemma~\ref{when there are independent coordinates with the same type over a}.
It will be used in the next section where we show that we can choose the coordinates to
be imaginaries defined by $\es$-definable equivalence relations on $M$ (rather than on $M^n$ for some $n > 1$),
and that dividing is controlled by these equivalence relations.

{\bf \em For the rest of this section we fix (an arbitrary) $0 < s \leq r$.}

\begin{rem}\label{remark on the first technical lemma}{\rm
(The intuition behind Lemma~\ref{different coordinates have different types over acl(emptyset)}.)
Let $\mcC$ be the structure where $C = \mbbN$ and the vocabulary of $\mcC$ is empty.
Let $G$ be the set of all 2-element subsets of $C$.
Turn $G$ into a graph $\mcG$ by saying that $a, b \in G$ are adjacent if and only if their intersection is a singleton.
Since $\mcC$ is $\omega$-categeorical and stable and $\mcG$ is interpretable in $\mcC$ (without parameters) it follows
that $\mcG$ is $\omega$-categorical and stable\footnote{
For example, by \cite[Theorem~7.3.8]{Hod} and \cite[Ch.~III, Lemma~6.7]{She}}, 
in fact superstable with SU-rank 2.\footnote{
Which follows by a straightforward argument using the definition of dividing.}
However, $\mcG$ is not homogeneous, because it is easy to see that the following two triples of elements
from $G$ satisfy the same quantifier-free formulas, but not the same formulas with quantifiers:
$(\{1, 2\}, \{2, 3\}, \{1, 3\}), (\{1, 2\}, \{1, 3\}, \{1, 4\})$.
Note that the intersection of the elements in the first triple is empty, but the intersection of the elements in the second triple
is nonempty.

The idea of the proof of Lemma~\ref{different coordinates have different types over acl(emptyset)}
is as follows, where we let $\overline{\crd}_s(a)$ abbreviate
`$\crd_s(a) \setminus C_{s-1}$' : If $a \in M$ and $c_1, c_2 \in \overline{\crd}_s(a)$ satisfy the
premisses of the lemma, and $E_{s-1}(c_1, c_2)$, then we can
find $a, a', a'', a^* \in M$ such that
$aa'' \equiv aa^*$ and $a'a'' \equiv a'a^*$, but $aa'a'' \not\equiv aa'a^*$.
This is done by choosing the elements in such a way that
$\overline{\crd}_s(a) \cap \overline{\crd}_s(a') \cap \overline{\crd}_s(a'') = \es$
and $\overline{\crd}_s(a) \cap \overline{\crd}_s(a') \cap \overline{\crd}_s(a^*) \neq \es$.

The proof of Lemma~3.9 in \cite{Kop-one-based} builds on the same idea.
But in its context, $s = 1$ so all elements of $C_s$ have SU-rank 1. 
Then, by \cite[Theorem~5.1]{AK}, 
the ``canonically embedded structure'' (in $\mcM\meq$) with universe $C_1$, is,
modulo ``dividing out by the relation $\acl(x) = \acl(y)$'', a 
reduct of a binary random structure.
This simplified the arguments in the proof of \cite[Lemma~3.9]{Kop-one-based}.
Here we use only (besides the given coordinatization) properties of forking/dividing and, in particular, 
the independence theorem for simple structures; but the arguments become more intricate.
}\end{rem}

\begin{lem}\label{different coordinates have different types over acl(emptyset)}
If $a \in M$, $c_1, c_2 \in \crd_s(a) \setminus C_{s-1}$, $ac_1 \equiv ac_2$ and $c_1 \underset{\crd_{s-1}(c_1)}{\ind} c_2$, then
\[
\tp\big(c_1 / \acl(\crd_{s-1}(c_1))\big) 
\neq \tp\big(c_2 / \acl(\crd_{s-1}(c_1))\big), \ \text{ and hence}
\]
\[
\tp\big(c_1 / \acl(\crd_{s-1}(a))\big) 
\neq \tp\big(c_2 / \acl(\crd_{s-1}(a))\big).
\]
\end{lem}

\noindent
{\bf Proof.}
For a contradiction suppose that there are $a \in M$ and $c_1, c_2 \in \crd_s(a) \setminus C_{s-1}$ such that
\begin{equation}\label{ac-1 equiv ac-2}
ac_1 \equiv ac_2,  \quad \ c_1 \underset{\crd_{s-1}(c_1)}{\ind} c_2, \ \text{ and}
\end{equation}
\begin{equation}\label{c-1 / acl(empty) = c-2 / acl(empty)}
\tp(c_1 / \acl(\crd_{s-1}(c_1))) = \tp(c_2 / \acl(\crd_{s-1}(c_1))).
\end{equation}
Note that this implies that $\crd_{s-1}(c_1) = \crd_{s-1}(c_2)$, so
\begin{equation}\label{E(c-1, c-2)}
E_{s-1}(c_1, c_2).
\end{equation}
By~(\ref{ac-1 equiv ac-2}) there is $c^*_1 \in C$ such that 
\begin{equation}\label{ac-1c-2 equiv ac-2-c*-1}
ac_1c_2 \equiv ac_2c^*_1.
\end{equation}
Then 
\begin{equation}\label{c-1 ind c-2 etc}
c^*_1 \in \crd_s(a) \setminus C_{s-1} \ \ \text{ and } \ \ c_2 \underset{\crd_{s-1}(c_1)}{\ind} c^*_1.
\end{equation}
From~(\ref{E(c-1, c-2)}),~(\ref{ac-1c-2 equiv ac-2-c*-1}) and Fact~\ref{fact about E-s} we also get
\begin{align}\label{c-2 / acl(empty) = c-1-* / acl(empty)}
E_{s-1}(c_2, c^*_1).
\end{align}
By~(\ref{ac-1 equiv ac-2}),~(\ref{c-1 ind c-2 etc}),~(\ref{c-2 / acl(empty) = c-1-* / acl(empty)}) and the independence theorem
there is $c'_2 \in C_s \setminus C_{s-1}$ such that
\begin{equation}\label{c-1c'-2 equiv c-1c-2 etc}
c_1c'_2 \equiv c_1c_2 \equiv c_2c'_2 \ \text{ and } \ 
c'_2 \underset{\crd_{s-1}(c_1)}{\ind} c_1, c_2.
\end{equation}
In addition, we may, without loss of generality, assume that 
\begin{equation}\label{c'-2 ind from a over c-1c-2}
c'_2 \underset{c_1c_2}{\ind} a,
\end{equation} 
because if this is not the case
then we can replace $c'_2$ by a realization of a nondividing extension of 
$\tp(c'_2 / c_1, c_2)$ to $\{a, c_1, c_2\}$
(and recall that $\crd_{s-1}(c_1) \subseteq \acl(c_1)$).

Since (by Fact~\ref{existence of coordinates}~(iii)) $\su(c'_2 / \crd_{s-1}(c_1)) \geq 1$, it follows 
from $c'_2 \underset{\crd_{s-1}(c_1)}{\ind} c_1c_2$ (see~(\ref{c-1c'-2 equiv c-1c-2 etc})) that
$c'_2 \notin \acl(c_1, c_2)$.
From this together with~(\ref{c'-2 ind from a over c-1c-2}) we get
\begin{equation}\label{c'-2 is not in crd(a)}
c'_2 \notin \acl(a) \ \text{ so } \ c'_2 \notin \crd_s(a).
\end{equation}
From~(\ref{c-1c'-2 equiv c-1c-2 etc}),~(\ref{c'-2 ind from a over c-1c-2}) and transitivity, we get
\begin{equation}\label{c'-2 ind from a over}
c'_2 \underset{\crd_{s-1}(c_1)}{\ind} a.
\end{equation}
By~(\ref{c-1c'-2 equiv c-1c-2 etc}) there are $a', a'' \in M$ such that
\begin{equation}\label{a'c-1c'-2 equiv ac-1c-2 etc}
a'c_1c'_2 \equiv ac_1c_2 \equiv a''c_2c'_2.
\end{equation}
By considering nondividing extensions if necessary we may assume, without loss of generality, that
\begin{equation}\label{c-1 c-2 ind a etc}
a' \underset{c_1c'_2}{\ind} a \ \ \text{ and } \ \ a'' \underset{c_2c'_2}{\ind} aa'.
\end{equation}
Before continuing, observe that for every $c \in C_s$, $\crd(c) = \crd_s(c) \subseteq C_s$, because of Fact~\ref{existence of coordinates}~(iii).

\begin{claim}\label{first claim about coordinates}
\begin{align}
&\crd(a) \cap \crd(a')  \ = \ \crd(c_1), \label{first equality}\\
&\crd(a) \cap \crd(a'')  \ = \ \crd(c_2), \  \text{ and}
\label{second equality}\\
&\crd(a') \cap \crd(a'') \ = \ \crd(c'_2). \label{third equality}
\end{align}
\end{claim}

\noindent
{\bf Proof of the claim.}
First note that by the choice of $a, c_1$ and $c_2$, and by~(\ref{a'c-1c'-2 equiv ac-1c-2 etc}),
we get $c_1 \in \crd(a) \cap \crd(a') \cap C_s$. 
Hence $\crd(c_1) \subseteq \crd(a) \cap \crd(a')$.
From~(\ref{c-1 c-2 ind a etc}) and Lemma~\ref{a more general independence criterium} we get
\begin{equation}\label{intersection in acl(c-1 c'-2)}
\crd(a) \cap \crd(a') \ \subseteq \ \acl(c_1, c'_2).
\end{equation}
Regarding~(\ref{first equality}), it remains to prove that $\crd(a) \cap \crd(a') \ \subseteq \ \crd(c_1)$.
Suppose that $d \in \crd(a) \cap \crd(a')$.
By~(\ref{intersection in acl(c-1 c'-2)}) and 
Fact~\ref{fact about C and crd}~(i),
$d \in \acl(c_1)$ or $d \in \acl(c'_2)$.
If $d \in \acl(c_1)$ then we have $d \in \crd(c_1)$.

Suppose that $d \in \acl(c'_2)$.
Hence $d \in \crd(a)  \cap \crd(c'_2)$.
From~(\ref{c'-2 ind from a over}) we have 
$c'_2 \underset{\crd_{s-1}(c_1)}{\ind} a$, so by
Lemma~\ref{a more general independence criterium}
we get $d \in \acl(\crd_{s-1}(c_1))$
and hence (by the definition of $\crd_{s-1}$) $d \in \crd(c_1)$.
Thus we have proved~(\ref{first equality}).

Observe that~(\ref{first equality}) and Lemma~\ref{a more general independence criterium} imply that 
\begin{equation}\label{a ind from a' over c-1}
a \underset{c_1}{\ind} a'.
\end{equation}
If $c_2 \in \acl(a')$ then, as $c_2 \in \acl(a)$, it follows from~(\ref{a ind from a' over c-1}) and
Lemma~\ref{a more general independence criterium}
that $c_2 \in \acl(c_1)$, but this contradicts~(\ref{ac-1 equiv ac-2}).
Hence,
\begin{equation}\label{c-2 is not in acl(a')}
c_2 \notin \acl(a')
\end{equation}

Now we prove~(\ref{second equality}).
From~(\ref{a'c-1c'-2 equiv ac-1c-2 etc}) it follows that 
$c_2 \in \crd(a) \cap \crd(a'')$, so 
$\crd(c_2) \subseteq \crd(a) \cap \crd(b)$.
It remains to prove that if $d \in \crd(a) \cap \crd(b)$ then $d \in \acl(c_2)$.
So suppose that $d  \in \crd(a) \cap \crd(b)$.
By~(\ref{c-1 c-2 ind a etc})
and Lemma~\ref{a more general independence criterium}, 
$d \in \acl(c_2, c'_2)$, so by Fact~\ref{fact about C and crd}~(i),
$d \in \acl(c_2)$ or $d \in \acl(c'_2)$.
If $d \in \acl(c_2)$ then we are done, so suppose that $d \in \acl(c'_2)$.

First assume that $d \in C_s \setminus C_{s-1}$.
Recall that, by Fact~\ref{fact about C and crd}~(ii), $(C_s \setminus C_{s-1}, \cl)$, where `$\cl$' is `$\acl$' restricted to 
$C_s \setminus C_{s-1}$, is a trival pregeometry.
By assumption, $d \in \acl(c'_2)$, so
(by the ``exchange property'' of pregeometries) 
$c'_2 \in \acl(d)$ and hence $c'_2 \in \acl(a)$, contradicting~(\ref{c'-2 is not in crd(a)}).

Hence we must have $d \in C_{s-1}$.
By assumption we have $d \in \crd(c'_2) \cap \crd(a)$.
This together with~(\ref{c'-2 ind from a over}) and Lemma~\ref{a more general independence criterium}
implies that $d \in \crd_{s-1}(c_1)$.
By~(\ref{E(c-1, c-2)}), $E_{s-1}(c_1, c_2)$, so $\crd_{s-1}(c_1) = \crd_{s-1}(c_2)$ and therefore $d \in \crd_{s-1}(c_2)$.
Thus~(\ref{second equality}) is proved.

It remains to prove~(\ref{third equality}).
By~(\ref{a'c-1c'-2 equiv ac-1c-2 etc}), $c'_2 \in \crd(a') \cap \crd(a'')$,
so $\crd(c'_2) \subseteq \crd(a') \cap \crd(a'')$.
It remains to prove that
if $d \in \crd(a') \cap \crd(a'')$ then $d \in \crd(c'_2)$.
Suppose that $d \in \crd(a') \cap \crd(a'')$.
Then, from~(\ref{c-1 c-2 ind a etc}) and Lemma~\ref{a more general independence criterium},
we get $d \in \acl(c_2, c'_2)$.
By Fact~\ref{fact about C and crd}~(i), $d \in \acl(c_2)$ or $d \in \acl(c'_2)$.
If $d \in \acl(c'_2)$ then we are done, so suppose that $d \in \acl(c_2)$.

First assume that $d \in C_s \setminus C_{s-1}$.
As $C_s \setminus C_{s-1}$ is a trivial pregeometry (with `$\acl$' restricted to  $C_s \setminus C_{s-1}$)
and $d \in \acl(c_2)$ 
we  get $c_2 \in \acl(d) \subseteq \acl(a')$, which contradicts~(\ref{c-2 is not in acl(a')}).

Hence we have $d \in C_{s-1}$. Then $d \in \crd_{s-1}(c_2)$.
By~(\ref{E(c-1, c-2)}), $E_{s-1}(c_1, c_2)$
and by~(\ref{c-1c'-2 equiv c-1c-2 etc}) we get $E_{s-1}(c_2, c'_2)$, so $\crd_{s-1}(c_2) = \crd_{s-1}(c'_2)$.
Therefore $d \in \crd_{s-1}(c'_2)$.
This concludes the proof of Claim~\ref{first claim about coordinates}.
\hfill $\square$
\\

\noindent
By~(\ref{ac-1 equiv ac-2}) there is $d \in M$ such that
\begin{equation}\label{ac-1d equiv ac-2a''}
ac_1d \equiv ac_2a'',
\end{equation}
so in particular, $c_1 \in \crd_s(d) \setminus C_{s-1}$.
By~(\ref{ac-1 equiv ac-2}) and~(\ref{a'c-1c'-2 equiv ac-1c-2 etc}) we have
$a'c_1 \equiv a'c'_2$, so there is $e \in M$ such that
\begin{equation}\label{a'b-1e equiv a'c'-2a''}
a'c_1e \equiv a'c'_2a'',
\end{equation}
so in particular, $c_1 \in \big(\crd_s(d) \cap \crd_s(e)\big) \setminus C_{s-1}$.
By~(\ref{a'b-1e equiv a'c'-2a''}),
~(\ref{a'c-1c'-2 equiv ac-1c-2 etc}),
~(\ref{ac-1 equiv ac-2}),
~(\ref{a'c-1c'-2 equiv ac-1c-2 etc})
and~(\ref{ac-1d equiv ac-2a''}),
in the mentioned order, we have
\[
c_1e \ \equiv \ c'_2a'' \ \equiv \ c_2a \ \equiv \ c_1a \ \equiv \ c_2a'' \ \equiv \ c_1d.
\]
Hence $c_1e \equiv c_1d$ and by Lemma~\ref{if c belongs to crd(a) and crd(b) then a and b have the same type over acl(c)} 
we get
\begin{equation}\label{d / acl(c-1) = e / acl(c-1)}
\tp(d / \acl(c_1)) = \tp(e / \acl(c_1)).
\end{equation}
From~(\ref{second equality}),~(\ref{third equality}) and
Lemma~\ref{a more general independence criterium}
we get
\[
a \underset{c_2}{\ind} a'' \ \ \text{ and } \ \ a' \underset{c'_2}{\ind} a'',
\]
which together with~(\ref{ac-1d equiv ac-2a''}) and~(\ref{a'b-1e equiv a'c'-2a''}) gives
\begin{equation}\label{a ind from d over c-1 etc}
a \underset{c_1}{\ind} d \ \ \text{ and } \ \ a' \underset{c_1}{\ind} e.
\end{equation}
By (\ref{a ind from a' over c-1}),~(\ref{d / acl(c-1) = e / acl(c-1)}), (\ref{a ind from d over c-1 etc}) 
and the independence theorem there is $a^* \in M$ such that
\begin{equation}\label{ac-1a* equiv ac-2a'' etc}
ac_1a^* \equiv ac_2a'' \ \ \text{ and } \ \ a'c_1a^* \equiv a'c'_2a''.
\end{equation}
This together 
with~(\ref{second equality}) and~(\ref{third equality}) 
implies that
\begin{align}\label{crd(a) intersected with crd(a*) is c-1 etc}
&\crd(a) \cap \crd(a^*) = \crd(c_1) \ \ \text{ and} \\
&\crd(a') \cap \crd(a^*) = \crd(c_1). \nonumber
\end{align}
Hence
\begin{equation}\label{c-1 is in crd(a), crd(a') and crd(a*)}
c_1 \in \crd(a) \cap \crd(a') \cap \crd(a^*) \cap (C_s \setminus C_{s-1}).
\end{equation}
By~(\ref{ac-1 equiv ac-2})
and~(\ref{c-1c'-2 equiv c-1c-2 etc}), $\{c_1, c_2, c'_2\}$ is an independent set over 
$\crd_{s-1}(c_1)$. 
Hence $\acl(c_1) \cap \acl(c_2) \cap \acl(c'_2) \cap (C_s \setminus C_{s-1}) = \es$.
Now Claim~\ref{first claim about coordinates}
implies that
\begin{equation}\label{a, a', a'' has empty intersection of coordinates at level s}
\crd(a) \cap \crd(a') \cap \crd(a'') \cap (C_s \setminus C_{s-1}) = \es.
\end{equation}
Since $a, a', a'', a^* \in M$ and $\mcM$ is a binary structure with elimination of quantifiers, it follows 
from~(\ref{ac-1a* equiv ac-2a'' etc}) that
\begin{equation}\label{aa'a'' equiv aa'a*}
aa'a'' \equiv aa'a^*.
\end{equation}
But this contradicts~(\ref{c-1 is in crd(a), crd(a') and crd(a*)}) and~(\ref{a, a', a'' has empty intersection of coordinates at level s}),
because the relation
``$\crd(x) \cap \crd(y) \cap \crd(z) \cap (C_s \setminus C_{s-1})$ is nonempty''
is $\es$-definable in $\mcM$.
This concludes the proof of Lemma~\ref{different coordinates have different types over acl(emptyset)}.
\hfill $\square$
\\

\noindent
Before proving our next main lemma we need the following auxilliary lemma:

\begin{lem}\label{about coordinates and equivalence classes}
Let $a \in M$, $c \in \crd_s(a) \setminus C_{s-1}$ and $p(a, c) = \tp(a, c)$. \\
(i) Suppose that $\acl(c) \cap p(a, \mcM\meq)$ has nonempty intersection with
each one of the $E_{s-1}$-equivalence classes $X_1, \ldots, X_n$. 
Furthermore, suppose that $a' \in M$ and $E_{s-1}(a, a')$ (so in particular $\crd_{s-1}(a) = \crd_{s-1}(a')$).
Then there is $c' \in \crd_s(a')$ such that $p(a', c')$ and 
$\acl(c') \cap p(a', \mcM\meq)$ has nonempty intersection with 
all of $X_1, \ldots, X_n$.\\
(ii) Suppose that $a' \in M$, $c' \in \crd_s(a') \setminus C_{s-1}$ and $p(a', c')$. Then 
$\acl(c') \cap p(a', \mcM\meq)$ 
has nonempty intersection with the same number of $E_{s-1}$-equivalence classes 
as $\acl(c) \cap p(a, \mcM\meq)$ has. \\
(iii) Suppose that $X_1, \ldots, X_n$ is an enumeration of all $E_{s-1}$-equivalence classes with which 
$\acl(c) \cap p(a, \mcM\meq)$ has nonempty intersection. 
Furthermore suppose that $a' \in M$ and $E_{s-1}(a, a')$.
If $c' \in \crd_s(a') \setminus C_{s-1}$, $p(a', c')$ and $E_{s-1}(c, c')$, then 
$\acl(c') \cap p(a', \mcM\meq)$ 
has nonempty intersection with all of $X_1, \ldots, X_n$.
\end{lem}

\noindent
{\bf Proof.}
Let $a \in M$, $c \in \crd_s(a) \setminus C_{s-1}$ and $p(a, c) = \tp(a, c)$.
In this proof we abbreviate $E_{s-1}$ by $E$.

(i) We first note that $c_E$ may, strictly speaking, be an element of $(M\meq)\meq$.
But since $\mcM\meq$ has elimination of imaginaries we may identify $c_E$ with an element of $M\meq$.
By slight abuse of terminology, we also denote the sort of $c_E$ by $E$.
Let $\acl(c) \cap p(a, \mcM\meq) = \{c_1, \ldots, c_n\}$ and, for each $i = 1, \ldots, n$, let
$X_i = [c_i]_E$.
From the definition of $E (= E_{s-1})$ it follows that $(c_i)_E \in \acl(\crd_{s-1}(a))$ for all $i = 1, \ldots, n$.
Let $\varphi(x, z_1, \ldots, z_n)$ be a formula in the language of $\mcM\meq$ which expresses the following condition:
\begin{itemize}
\item[] ``each one of $z_1, \ldots, z_n$ is of sort $E$ and \\
$\exists y \Big( p(x,y) \ \wedge \ 
\forall u \Big( \big( p(x,u) \ \wedge \ u \in \acl(y) \big) \ \rightarrow$ \\ 
$\text{for some $1 \leq i \leq n$, $u$ belongs to the $E$-class represented by $z_i$} \Big)\Big)$''.
\end{itemize}
Then $\mcM\meq \models \varphi(a, (c_1)_E, \ldots, (c_n)_E)$. 
Let $a' \in M$ be such that $E(a, a')$.
Then $(c_i)_E \in \acl(\crd_{s-1}(a)) = \acl(\crd_{s-1}(a'))$ for all $i$, and
\[
\tp\big(a / \acl(\crd_{s-1}(a))\big) \ = \ \tp\big(a' / \acl(\crd_{s-1}(a'))\big).
\]
Hence we get
$\mcM\meq \models \varphi(a', (c_1)_E, \ldots, (c_n)_E)$. 
Thus there is $c' \in \crd_s(a')$ such that $\mcM\models p(a', c')$ and 
$\acl(c') \cap p(a, \mcM\meq)$ 
has nonempty intersection with $X_i$
for each $i = 1, \ldots, n$.

(ii) The assumption that $p(a, c)$ and $p(a', c')$ gives $ac \equiv a'c'$ so
there is an automorphism of $\mcM\meq$ which takes $ac$ to $a'c'$.
The conclusion follows from this.

(iii) Let $X_1, \ldots, X_n$ be an enumeration of all $E$-classes with which
$\acl(c) \cap p(a, \mcM\meq)$ has nonempty intersection.
Suppose that $a' \in M$, $E(a, a')$, 
$c' \in \crd_s(a')$, $p(a', c')$ and $E(c, c')$.
Using part~(ii) we can enumerate all $E$-classes with which
$\acl(c') \cap p(a', \mcM\meq)$ has nonempty intersection as $X'_1, \ldots, X'_n$.
Without loss of generality, assume that $X_1 = X'_1$ and $c, c' \in X_1$.
By part~(i), there is $c'' \in \crd_s(a') \setminus C_{s-1}$ such that $p(a', c'')$ and 
$\acl(c'') \cap p(a', \mcM\meq)$ 
has nonempty intersection with all $X_1, \ldots, X_n$. 
In particular, $\acl(c'') \cap p(a', \mcM\meq)$ has nonempty intersection with $X_1$.
Let $c^* \in \acl(c'') \cap p(a', \mcM\meq) \cap X_1$ (so in particular $c^* \in C_s \setminus C_{s-1}$).
As, by Fact~\ref{fact about C and crd}, 
$C_s \setminus C_{s-1}$ is a trivial pregeometry, with `$\acl$' restricted to $C_s \setminus C_{s-1}$,
we get $\acl(c^*) = \acl(c'')$.
Consequently $\acl(c^*) \cap p(a', \mcM\meq)$ has nonempty intersection with all $X_1, \ldots, X_n$.
By the choice of $c^*$ we have $a'c^* \equiv a'c'$ and $E(c^*, c')$.
Hence Lemma~\ref{different coordinates have different types over acl(emptyset)}
implies that $c^* \underset{\crd_{s-1}(c')}{\nind} c'$.
Since, by Fact~\ref{rank 1 over its coordinates on the previous level},
$\su(c' / \crd_{s-1}(c')) = 1$, we get $c' \in \acl(\{c^*\} \cup \crd_{s-1}(c'))$.
By Fact~\ref{fact about C and crd}~(i), we get
$c' \in \acl(c^*)$ or $c' \in \acl(\crd_{s-1}(c'))$.
But as $\su(c' / \crd_{s-1}(c')) = 1$ we must have $c' \in \acl(c^*)$.
Since $C_s \setminus C_{s-1}$ is a trivial pregeometry we get $\acl(c') = \acl(c^*)$.
Then $\acl(c') \cap p(a', \mcM\meq)$ has nonempty intersection with all
$X'_1, \ldots, X'_n, X_1, \ldots, X_n$, which, by part~(ii)
and the choice of $X_1, \ldots, X_n$ and $X'_1, \ldots, X'_n$, implies that $\{X_1, \ldots, X_n\} = \{X'_1, \ldots, X'_n\}$.
\hfill $\square$

\begin{rem}\label{remark on the second technical lemma}{\rm
(The intuition behind Lemma~\ref{when there are independent coordinates with the same type over a}.)
Let $\mcC = (\mbbN, E)$, where $E$ is interpreted as an equivalence relation with two infinite equivalence classes. 
Let us assume that one of the classes contains all even numbers and the other all odd numbers.
Let 
\[
G = \{\{n, m\} : n \in \mbbN \text{ is even and } m \in \mbbN \text{ is odd}\}.
\]
Turn $G$ into a graph $\mcG$ by letting $a, b \in G$ be adjacent if and only if their intersection is a singleton.
Since $\mcC$ is $\omega$-categorical and stable, and $\mcG$ is interpretable in $\mcC$ (without parameters)
it follows that $\mcG$ is $\omega$-categorical and stable, in fact superstable of SU-rank~2.
Without going into the details, we may assume, without loss of generality,
that $C \ (= \mbbN)$ is a $\es$-definable subset of $G\meq$ and that
the equivalence relation $E$ on $C$ is $\es$-definable in $\mcG\meq$.
Consider the following two quadruples of elements from $G$:
\[
(\{1, 2\}, \{1, 4\}, \{3, 6\}, \{3, 8\}), \ (\{1, 2\}, \{1, 4\}, \{3, 6\}, \{5, 6\}).
\]
Clearly, the two quadruples satisfy the same quantifier-free formulas.
Note that $\{1, 2\}$ and $\{1, 4\}$ have a common element in the $E$-class of odd numbers, and
the same is true for $\{3, 6\}$ and $\{3, 8\}$. 
Hence the first quadruple above satisfies the formula $\varphi(x_1, x_2, x_3, x_4)$ which expresses
``there are $u, v \in C$ such that $E(u, v)$, $x_1 \cap x_2 = \{u\}$ and $x_3 \cap x_4 = \{v\}$''.
But the second quadruple does not satisfy this formula. 
Since all elements in the two quadruples above are ``real'' elements of $G\meq$ (i.e. belong to $G$), it follows
that there is a formula in the (graph) language of $\mcG$ which is satisfied by the first quadruple, but not by the second.
Thus $\mcG$ is not homogeneous.

The idea of the proof of Lemma~\ref{when there are independent coordinates with the same type over a}
is the following:
If $a \in M$, $c_1, c_2 \in \crd_s(a) \setminus C_{s-1}$, 
$c_1 \underset{\crd_{s-1}(a)}{\ind} c_2$, and $ac_1 \equiv ac_2$,
then we can  find $a^*, b^*, a', b', b'' \in M$ such that
\[
a^*b^*b' \equiv a^*b^*b'' \ \text{ and } \ a'b' \equiv a'b'', \ \text{ but } \
a^*b^*a'b' \not\equiv a^*b^*a'b''.
\]
This is done by choosing the elements so that, with $p = \tp(a, c_1)$,
there are $c, d \in C_s \setminus C_{s-1}$ such that 
$E_{s-1}(c, d)$, $p(a^*, c)$, $p(b^*, c)$, $p(a', d)$ and $p(b', d)$, but no such $c$ and $d$ exist
if we replace $b'$ by $b''$.
In finding such elements we use 
Lemma~\ref{different coordinates have different types over acl(emptyset)},
which implies that
$\neg E_{s-1}(c_1, c_2)$, where `$E_{s-1}$' plays the role of  `$E$' in  $\mcG\meq$.

The same idea is behind the proof of \cite[Proposition~4.4]{Kop16b}, as becomes
apparent in the last page of that proof. However, in the context of \cite{Kop16b} one can assume
that $s = 1$, and then all $c \in C_s$ have SU-rank 1. Moreover, one can assume (in \cite{Kop16b})
that for all $c, d \in C_s$, if $d \in \acl(c)$, then $c = d$, and that the ``canonically embedded'' structure (in $\mcM\meq$)
with universe $C_s$ is a binary random structure (by \cite[Theorem~5.1]{AK}
and some additional observations in \cite[Fact~3.6]{Kop16b}).
In the present context, the arguments in the more specialized situation of \cite{Kop16b} 
are replaced by dividing/forking arguments.
}\end{rem}

\begin{lem}\label{when there are independent coordinates with the same type over a}
For all $a \in M$ and all $c_1, c_2 \in \crd_s(a) \setminus C_{s-1}$, 
if $c_1 \underset{\crd_{s-1}(a)}{\ind} c_2$ then $ac_1 \not\equiv ac_2$.
\end{lem}

\noindent
{\bf Proof.}
Towards a contradiction suppose that there are $a \in M$ and 
$c_1, c_2 \in \crd_s(a) \setminus C_{s-1}$ such that $c_1  \underset{\crd_{s-1}(a)}{\ind} c_2$ and $ac_1 \equiv ac_2$.
Let 
\[
q(x) = \tp(a) \ \ \text{ and } \ \ p(x, y) = \tp(a, c_1).
\]
Note that if $p(a', c)$ then $c \in \crd_s(a') \setminus C_{s-1}$.
So for every  $a' \in M$ which realizes $q$ there are $c, c' \in \crd_s(a') \setminus C_{s-1}$ 
such that $c  \underset{\crd_{s-1}(a')}{\ind} c'$ and both $a'c$ and $a'c'$ realize $p$.
Also, for all $a'$ and $c$ such that $a'c$ realizes $p$ there is $c'$ such that $a'c'$ realizes $p$ and
$c  \underset{\crd_{s-1}(a')}{\ind} c'$.

Choose any  $c \in \crd_s(a) \setminus C_{s-1}$ such that $ac$ realizes $p$.
Let $b \in M$ realize a nondividing extension of $\tp\big(a / \acl\big(\{c\} \cup \crd_{s-1}(a)\big)\big)$ 
to $\{a\} \cup \acl\big(\{c\} \cup \crd_{s-1}(a)\big)$.
Then
\begin{equation}\label{c in the intersection of crd(a) and crd(b)}
a \underset{\substack{\{c\} \cup \\ \crd_{s-1}(a)}}{\ind} b, \quad E_{s-1}(a, b) \quad \text{ and } \quad  p(a, c) \wedge p(b, c).
\end{equation}
By the choice of $p$,~(\ref{c in the intersection of crd(a) and crd(b)}) and
Lemma~\ref{a more general independence criterium} 
we get
\[
p(a, \mcM\meq) \cap p(b, \mcM\meq) \ \subseteq \ \crd(a) \cap \crd(b) \ \subseteq \ \acl\big(\{c\} \cup \crd_{s-1}(a)\big).
\]
Let $d \in p(a, \mcM\meq) \cap p(b, \mcM\meq)$.
By Fact~\ref{fact about C and crd}~(i), $d \in \acl(c)$ or $d \in \acl(\crd_{s-1}(a))$.
In the later case $d \in C_{s-1}$, because of Fact~\ref{existence of coordinates}~(iii), and
this contradicts that $p(a, \mcM\meq) \subseteq C_s \setminus C_{s-1}$.
Hence $d \in \acl(c)$, so we have proved that
\begin{equation}\label{p(a, M) cap p(b, M) is a subset of acl(c)}
c \in p(a, \mcM\meq) \cap p(b, \mcM\meq) \ \subseteq \ \crd(c).
\end{equation}
Let $a' \in M$ realize a nondividing extension of $\tp\big(a / \acl\big(\crd_{s-1}(a)\big)\big)$ 
to $\{a, b\} \cup \acl\big(\crd_{s-1}(a)\big)$.
Then 
\begin{equation}\label{E(a, a') in second main lemma}
E_{s-1}(a, a'), \quad \quad a' \underset{\crd_{s-1}(a)}{\ind} ab,
\end{equation}
and by Lemma~\ref{a more general independence criterium}
and Fact~\ref{fact about C and crd}~(iii),
\begin{equation}\label{crd(a) has empty intersection with crd(ab)}
\crd(a') \cap (\crd(a) \cup \crd(b)) = \crd_{s-1}(a).
\end{equation}
By Lemma~\ref{about coordinates and equivalence classes}~(i)
there is $c' \in \crd_s(a')  \setminus C_{s-1}$ such that $p(a', c')$ and $E_{s-1}(c, c')$.
As explained in the beginning of the proof, there is $c'' \in \crd_s(a) \setminus C_{s-1}$ 
such that $p(a', c'')$ and $c' \underset{\crd_{s-1}(a')}{\ind} c''$.
By~(\ref{E(a, a') in second main lemma}),
$\crd_{s-1}(a) = \crd_{s-1}(a')$ and therefore
\begin{equation}\label{c' ind from c'' over crd(a')}
c' \underset{\crd_{s-1}(a)}{\ind} c''.
\end{equation}
Let $b' \in M$ realize a nondividing extension of 
\[
\tp(a' / \{c'\} \cup \acl(\crd_{s-1}(a)))
\ \text{ to } \ \{a', a, b, c'\} \cup \acl(\crd_{s-1}(a)).
\]
Then
\begin{equation}\label{E(a', b') in second main lemma}
E_{s-1}(a', b'), \quad \quad a'ab\underset{\substack{ \{c'\} \cup \\ \crd_{s-1}(a)}}{\ind} b',
\end{equation}
and,
in the same way as we proved~(\ref{p(a, M) cap p(b, M) is a subset of acl(c)}), 
we get
\begin{equation}\label{c' in the intersection of crd(a') and crd(b')}
c' \in p(a', \mcM\meq) \cap p(b', \mcM\meq) \subseteq \crd(c').
\end{equation}
From~(\ref{crd(a) has empty intersection with crd(ab)}) 
and $c' \in \crd_s(a')$ we get $c' \underset{\crd_{s-1}(a)}{\ind} ab$, so 
by~(\ref{E(a', b') in second main lemma})
and transitivity of dividing we also have
\begin{equation}\label{b' is independent from ab}
ab \underset{\crd_{s-1}(a)}{\ind} b'.
\end{equation}
Since $p(a', c')$, $p(a', c'')$ and $E_{s-1}(a', b')$, there is $b'' \in M$ such that
\begin{equation}\label{a'c'b' equiv a'c''b''}
a'c'b' \equiv a'c''b'',   \ \text{ so } \ E_{s-1}(a', b'') \ \text{ and hence } \ E_{s-1}(b', b'').
\end{equation}
This together with~(\ref{E(a', b') in second main lemma}) implies that
\begin{equation}\label{a' and b'' independent over c''}
a' \underset{\substack{\{c''\} \cup \\ \crd_{s-1}(a)}}{\ind} b''.
\end{equation}
Note that since $E_{s-1}(a, b)$, $E_{s-1}(a, a')$, $E_{s-1}(a', b')$
and $E_{s-1}(b', b'')$, all the elements $a, a', b, b'$ and $b''$ have the same type over
$\acl(\crd_{s-1}(a))$.
By considering a nondividing extension of 
\[
\tp\big(b'' / \{a', c''\} \cup \acl(\crd_{s-1}(a))\big) \ \text{ to } \ \{a', c'', a, b, b'\} \cup \acl(\crd_{s-1}(a)),
\]
if necessary, we may, in addition, assume that
\begin{equation}\label{b'' ind from abb' over}
b'' \underset{\substack{\{a', c''\} \cup \\ \crd_{s-1}(a)}}{\ind} abb'.
\end{equation}
This together with~(\ref{a' and b'' independent over c''}) and transitivity gives
$b'' \underset{\substack{\{c''\} \cup \\ \crd_{s-1}(a)}}{\ind} abb'$.
By the choice of $c''$, $c'' \in \crd_s(a')$. Hence~(\ref{E(a, a') in second main lemma}) implies that
$c'' \underset{\crd_{s-1}(a)}{\ind} ab$, 
so by transitivity 
\begin{equation}\label{b'' ind from ab over}
b'' \underset{\crd_{s-1}(a)}{\ind} ab.
\end{equation}

\begin{claim}\label{claim in second main lemma}
$c'' \underset{\crd_{s-1}(a)}{\ind} b'$.
\end{claim}

\noindent
{\bf Proof of the claim.}
By~(\ref{E(a', b') in second main lemma}), 
Lemma~\ref{a more general independence criterium}
and facts~\ref{existence of coordinates}~(iii)
and~\ref{fact about C and crd}~(i),
\[
\big(\crd_s(a') \cap \crd_s(b')\big) \setminus C_{s-1} \subseteq \acl(c').
\]
Recall that we have chosen $c''$ so that $c' \underset{\crd_{s-1}(a)}{\ind} c''$.
Hence $c'' \notin \acl(c')$.
Since $c'' \in \crd_s(a')$ it follows that $c'' \notin \acl(b')$,
and hence
\[
c'' \notin \crd_s(b').
\]
Suppose, for a contradiction, that there is $d \in \big(\crd_s(c'') \cap \crd_s(b')\big) \setminus C_{s-1}$.
Since $C_s \setminus C_{s-1}$ is a trivial pregeometry (by Fact~\ref{fact about C and crd}~(ii)),
we get $c'' \in \acl(d)$, and hence $c'' \in \crd_s(b')$, contradicting what we obtained above.
It follows that $\crd(c'') \cap \crd(b') \subseteq C_{s-1}$, so
\[
\crd(c'') \cap \crd(b') = \crd_{s-1}(c'') \cap \crd_{s-1}(b').
\]
Since $a'$ and $b'$ have the same type over $\acl(\crd_{s-1}(a)) = \acl(\crd_{s-1}(a')) = \acl(\crd_{s-1}(b'))$,
it follows that $\crd_{s-1}(a') = \crd_{s-1}(b')$.
As $c'' \in \crd_s(a')$ 
we get $\crd_{s-1}(c'') \subseteq \crd_{s-1}(b')$.
Consequently, $\crd_{s-1}(c'') \cap \crd_{s-1}(b') = \crd_{s-1}(c'')$.
Since we proved that $\crd(c'') \cap \crd(b') = \crd_{s-1}(c'') \cap \crd_{s-1}(b')$
it follows from
Lemma~\ref{a more general independence criterium}
that $c'' \underset{\crd_{s-1}(c'')}{\ind} b'$
and hence $c'' \underset{\crd_{s-1}(a)}{\ind} b'$.
\hfill $\square$
\\

\noindent
On the line after~(\ref{b'' ind from abb' over})
we obtained
$b'' \underset{\substack{\{c''\} \cup \\ \crd_{s-1}(a)}}{\ind} abb'$, 
from which we get 
$b'' \underset{\substack{\{c''\} \cup \\ \crd_{s-1}(a)}}{\ind} b'$.
This together with Claim~\ref{claim in second main lemma}
and transitivity gives
\begin{equation}\label{b' and b'' are independent}
b' \underset{\crd_{s-1}(a)}{\ind} b''.
\end{equation}
We have $E_{s-1}(b', b'')$ and this implies that $\crd_{s-1}(b') = \crd_{s-1}(b'') = \crd_{s-1}(a)$ and
\[
\tp(b' / \acl(\crd_{s-1}(a))) \ = \ \tp(b'' / \acl(\crd_{s-1}(a))).
\]
It follows (from Fact~\ref{generalized omega-homogeneity}) that there are $a^+, b^+ \in M$ such that 
\begin{equation}\label{abb' and a+b+b'' have the same type}
\tp(a, b, b' / \acl(\crd_{s-1}(a))) \ = \ \tp(a^+, b^+, b'' / \acl(\crd_{s-1}(a))),
\end{equation}
which by~(\ref{b' is independent from ab}) implies that
\begin{equation}\label{a+b+ independent from b''}
a^+b^+ \underset{\crd_{s-1}(a)}{\ind} b''.
\end{equation}
By~(\ref{b' is independent from ab}),~(\ref{b' and b'' are independent}),~(\ref{abb' and a+b+b'' have the same type}),
~(\ref{a+b+ independent from b''})
and the independence theorem there are $a^*, b^* \in M$ such that
\begin{align}\label{existence of a*b*}
&\tp(a^*, b^*, b' /  \acl(\crd_{s-1}(a))) \ = \ \tp(a, b, b' / \acl(\crd_{s-1}(a))), \\
&\tp(a^*, b^*, b'' / \acl(\crd_{s-1}(a))) \ = \ \tp(a^+, b^+, b'' / \acl(\crd_{s-1}(a))), \ \text{ and} \nonumber \\
&a^*b^* \underset{\crd_{s-1}(a)}{\ind} b'b''. \nonumber
\end{align}
By considering a nondividing extension if necessary we may, in addition, assume that
\[
a^*b^* \underset{\crd_{s-1}(a)}{\ind} a'b'b''.
\]
From~(\ref{abb' and a+b+b'' have the same type}) 
and~(\ref{existence of a*b*}) we get
\[
a^*b^*b' \equiv a^*b^*b''.
\]
From~(\ref{a'c'b' equiv a'c''b''}) we have $a'b' \equiv a'b''$.
Since $a^*, b^*, a', b', b'' \in M$ where $\mcM$ is binary with elimination of quantifiers it follows that 
\begin{equation}\label{a*b*a'b' equiv a*b*a'b''}
a^*b^*a'b' \equiv a^*b^*a'b''.
\end{equation}
By~(\ref{abb' and a+b+b'' have the same type}) and~(\ref{existence of a*b*})
we have $\tp(a^*, b^* / \acl(\crd_{s-1}(a))) = \tp(a, b / \acl(\crd_{s-1}(a)))$.
Recall that $c \in p(a, \mcM\meq) \cap p(b, \mcM\meq)$.
Therefore
(and by Fact~\ref{fact about isolated types}~(iii))
 there is $c^* \in p(a^*, \mcM\meq) \cap p(b^*, \mcM\meq)$ such that $E_{s-1}(c, c^*)$.
We have chosen $c'$ so that, among other things,
$E_{s-1}(c, c')$ (see the line after~(\ref{crd(a) has empty intersection with crd(ab)})). 
As $E_{s-1}$ is an equivalence relation we get $E_{s-1}(c', c^*)$. 
These observations and~(\ref{c' in the intersection of crd(a') and crd(b')})
imply that the following statement, abbreviated $\varphi(x_1, x_2, x_3, x_4)$, is satisfied by $(a^*, b^*, a', b')$:
\begin{itemize}
\item[] ``There are $y_1, y_2 \in C_s \setminus C_{s-1}$ such that $E_{s-1}(y_1, y_2)$ and
$p(x_1, y_1)$, $p(x_2, y_1)$, $p(x_3, y_2)$ and $p(x_4, y_2)$.''
\end{itemize}
Note that $\varphi(x_1, x_2, x_3, x_4)$ can be expressed by a first-order formula in the language of $\mcM\meq$.
The next step is to show that $\varphi$ is not satisfied by $(a^*, b^*, a', b'')$.

Suppose that $d, e \in C_s \setminus C_{s-1}$ are such that
\[
p(a^*, d) \cap p(b^*, d) \ \text{ and } \ 
p(a', e) \cap p(b'', e).
\]
To prove that $\mcM\meq \not\models \varphi(a^*, b^*, a', b'')$ it suffices to show that $\neg E_{s-1}(d, e)$.
By the choice of $c^*$, (\ref{c in the intersection of crd(a) and crd(b)}) and~(\ref{existence of a*b*}),
we have $a^* \underset{\substack{\{c^*\} \cup \\ \crd_{s-1}(a)}}{\ind} b^*$ and therefore
\[
\crd(a^*) \cap \crd(b^*) \cap (C_s \setminus C_{s-1}) \subseteq \acl(c^*).
\]
Moreover, from (\ref{a' and b'' independent over c''}) it follows that
\[
\crd(a') \cap \crd(b'') \cap (C_s \setminus C_{s-1}) \subseteq \acl(c'').
\]
Therefore the assumptions about $d$ and $e$ imply that
\[
d \in \acl(c^*) \cap p(a^*, \mcM\meq) \ \text{ and } \ e \in \acl(c'') \cap p(a', \mcM\meq).
\]
Since $C_s \setminus C_{s-1}$ is a trivial pregeometry it follows that $c'' \in \acl(e)$, and hence $\acl(e) = \acl(c'')$.
Recall that $E_{s-1}(c', c^*)$. 
By Lemma~\ref{about coordinates and equivalence classes}~(iii),
there is $e' \in \acl(c') \cap p(a', \mcM\meq)$ such that $E_{s-1}(d, e')$.
By again using that $C_s \setminus C_{s-1}$ is a trivial pregeometry it follows
that $c' \in \acl(e')$, and consequently $\acl(c') = \acl(e')$.
Thus we have $\acl(e') = \acl(c')$ and $\acl(e) = \acl(c'')$, and 
by~(\ref{c' ind from c'' over crd(a')})
we have $c' \underset{\crd_{s-1}(a)}{\ind} c''$.
It follows that $e \underset{\crd_{s-1}(a)}{\ind} e'$.
By the choice of $e$ and $e'$ we also have $a'e \equiv a'e'$.
Therefore Lemma~\ref{different coordinates have different types over acl(emptyset)}
implies that $\neg E_{s-1}(e, e')$.
Since $E_{s-1}(d, e')$ we must have $\neg E_{s-1}(d, e)$.
Thus we have shown that 
\[
\mcM\meq \models \varphi(a^*, b^*, a', b'') \wedge \neg\varphi(a^*, b^*, a', b''),
\]
which contradicts~(\ref{a*b*a'b' equiv a*b*a'b''}).
This concludes the proof of Lemma~\ref{when there are independent coordinates with the same type over a}.
\hfill $\square$

\section{Coordinatization by equivalence relations}\label{Coordinatization by equivalence relations}

\noindent
Throughout this section we adopt Notation~\ref{remark about acl and dcl}.
Theorem~\ref{theorem coordinatizing M by equivalence classes}, below,
is slightly more general\footnote{Because we only assume that $\bar{c} \in C$ here.} 
than (a) --~(c) of the main results in the introduction.
Corollaries~\ref{information about nonbinary structures}
and~\ref{information about nonbinary primitive structures} have more general assumptions
than Theorem~\ref{theorem coordinatizing M by equivalence classes}
and are derived from its proof.

\begin{theor}\label{theorem coordinatizing M by equivalence classes}
Suppose that $\mcM$ is binary, simple, and homogeneous (hence supersimple with finite SU-rank).
Let $\mbR$ be the (finite) set of all $\es$-definable equivalence relations on~$M$.
\begin{itemize}
\item[(i)] For every $a \in M$, if $\su(a) = k$, then there are $R_1, \ldots, R_k \in \mbR$,
depending only on $\tp(a)$,  such that
$a \in \acl(a_{R_k})$, $\su(a_{R_1}) =~1$, 
$R_{i+1}$ refines $R_i$ and 
$\su(a_{R_{i+1}} / a_{R_i}) = 1$ for all $1 \leq i < k$
(or equivalently, $\su(a / a_{R_i}) = k-i$ for all $1 \leq i \leq k$).

\item[(ii)] Suppose that $a, b \in M$, $\bar{c} \in C$, and $a \underset{\bar{c}}{\nind} b$
(where we recall that $M \subseteq C \subseteq M\meq$).
Then there is $R \in \mbR$ such that $a \underset{\bar{c}}{\nind} a_R$ and $a_R \in \acl(b)$
(and hence $a_R \notin \acl(\bar{c})$).

\item[(iii)] Suppose that all binary $\es$-definable relations on $M$ are symmetric.
If $a, b \in M$, $\bar{c} \in C$, and $a \underset{\bar{c}}{\nind} b$,
then  there is $R \in \mbR$ such that $a \underset{\bar{c}}{\nind} a_R$ and $R(a, b)$ (and therfore $a_R \in \acl(b)$,
$a_R \notin \acl(\bar{c})$ and hence $\neg E(a, c)$ for every $c \in \bar{c}$).\footnote{
Note that the assumptions of part~(iii) imply that $Th(\mcM)$ has only one 1-type over $\es$.}
\end{itemize}
\end{theor}

\begin{rem}\label{remarks on the coordinatization by equivalence relations}{\rm
(i) Suppose that $a \in M$. The ``coordinatization by $R_1, \ldots, R_k$'' as in 
Theorem~\ref{theorem coordinatizing M by equivalence classes}~(i) may {\em not} be unique.
In other words, there may also be $\es$-definable equivalence relations $R'_1, \ldots, R'_k$ with the same
properties as $R_1, \ldots, R_k$ such that some $R'_i$ is (in a strong sense\footnote{
For example, it can happen, like with $\mcM$ in Section~\ref{Bipedes with bicoloured legs}, that
$R_i$ and $R'_i$ have only infinite classes but $R_i \cap R'_i$ has only singleton classes.}) not equivalent with $R_i$.
This is shown by the example $\mcM$ in Section~\ref{Bipedes with bicoloured legs}.\\
(ii) The conclusion in Theorem~\ref{theorem coordinatizing M by equivalence classes}~(ii)
{\em cannot} be strengthened so that it, in addition, says that $R(a, b)$.
This is also shown by the example $\mcM$ in Section~\ref{Bipedes with bicoloured legs}.
}\end{rem}

\noindent
The following two corollaries follow from an analysis of the proof
of Theorem~\ref{theorem coordinatizing M by equivalence classes},
which is given in Section~\ref{proofs of two corollaries}.

\begin{cor}\label{information about nonbinary structures}
Suppose that $\mcM$ is $\omega$-categorical, supersimple with finite SU-rank and with trivial dependence.
Also, suppose that part~(i) of Theorem~\ref{theorem coordinatizing M by equivalence classes}
does not hold for $\mcM$.
Then there are distinct $a_i, b_i \in M$, $i = 1, \ldots, 4$, 
such that $\tp(a_i, a_j) = \tp(b_i, b_j)$ for all $i, j$ and $\tp(a_1, \ldots, a_4) \neq \tp(b_1, \ldots, b_4)$.
So if $\mcM$ is, in addition, homogeneous, then it must have some relation symbol of arity 3 or 4.
\end{cor}

\begin{cor}\label{information about nonbinary primitive structures}
Suppose that $\mcM$ is $\omega$-categorical, supersimple with finite SU-rank and with trivial dependence.
Moreover, assume that $\mcM$ has no $\es$-definable equivalence relation on $M$ with infinitely many infinite equivalence classes.
If $\su(\mcM) > 1$ then there are distinct $a_i, b_i \in M$, $i = 1, \ldots, 4$, 
such that $\tp(a_i, a_j) = \tp(b_i, b_j)$, for all $i, j$, and $\tp(a_1, \ldots, a_4) \neq \tp(b_1, \ldots, b_4)$.
So if $\mcM$ is, in addition, homogeneous, then it must have some relation symbol of arity 3 or 4.
\end{cor}

\subsection{Proof of part~(i) of Theorem~\ref{theorem coordinatizing M by equivalence classes}}
\label{proof of part i}

In this subsection (and the next) we assume that $\mcM$ is binary, simple and homogeneous.
Moreover, we assume that $M \subseteq U \subseteq M\meq$, where $U$, $C$ and $C_i$, $i = 1, \ldots, h$,
are as in Assumption~\ref{assumption about C}.
Then we can use all results from sections~\ref{Preliminaries}
--~\ref{The main technical lemmas}.
The proof is carried out through a sequence of lemmas and is finished by the short argument after 
Lemma~\ref{the rank of an element over its predeccessor}.

\begin{lem}\label{being in the algebraic closure of an equivalence relation}
Suppose that $Q$ is a $\es$-definable equivalence relation on $M^n$.
Let $\bar{a} \in M^n$ and suppose that $b \in \acl(\bar{a}')$ for every $\bar{a}' \in [\bar{a}]_Q$.
Then $b \in \acl(\bar{a}_Q)$.
\end{lem}

\noindent
{\bf Proof.}
If $[\bar{a}]_Q$ is finite the $\acl(\bar{a}) = \acl(\bar{a}_Q)$ and the conclusion is immediate.
So suppose that $[\bar{a}]_Q$ is infinite.
For a contradiction suppose that $b \notin \acl(\bar{a}_Q)$.
The we find $\bar{a}'$ (in some elementary extension of $\mcM$) 
realizing a nonforking extension of $\tp(\bar{a} / \bar{a}_Q)$
to $\bar{a}_Q b$, so $\bar{a}' \underset{\bar{a}_Q}{\ind} b$.
As $\mcM\meq$ is $\omega$-saturated we may assume that $\bar{a}' \in M^n$.
Since $\tp(\bar{a}' / \bar{a}_Q) = \tp(\bar{a} / \bar{a}_Q)$ we have $\bar{a}' \in [\bar{a}]_Q$.
As  $\bar{a}' \underset{\bar{a}_Q}{\ind} b$ and $b \notin \acl(\bar{a}_Q)$, 
we get $b \notin \acl(\bar{a}')$, contradicting the assumption.
\hfill $\square$

\begin{defin}\label{definition of R-q}{\rm
Let $a \in M$, $c \in \crd(a)$, $q(u) = \tp(c)$ and $p(x, u) = \tp(a, c)$.
Define a relation on $M$ as follows:
\[
R_p(x, y) \ \Longleftrightarrow \  \Big( \neg q(x) \ \wedge \  \neg q(y) \Big) \ \vee \ 
\exists u, v \Big( p(x, u) \wedge p(x, u) \wedge \acl(u) = \acl(v) \Big).
\]
}\end{defin}

\begin{lem}\label{R-q is an equivalence relation}
The relation $R_p$, as in Definition~\ref{definition of R-q},
is an equivalence relation and is $\es$-definable.
\end{lem}

\noindent
{\bf Proof.}
By $\omega$-categoricity, $R_p$ is $\es$-definable.
It is straightforward to see that it is reflexive and symmetric, so it remains to show that it is transitive.
Suppose that $a, b, c \in M$, $R_p(a, b)$ and $R_p(b, c)$.
We assume that $a \neq b$, $a \neq c$, $b \neq c$, $q(a)$, $q(b)$, and $q(c)$, as the other cases are straightforward
and only use the definition of $R_p$.
By the definition of $R_p$, there are $i, j, k, l$  and $c_i, c_j, c_k, c_l$ such that
$p(a, c_i)$, $p(b, c_j)$, $p(b, c_k)$, $p(c, c_l)$,
$\acl(c_i) = \acl(c_j)$ and $\acl(c_k) = \acl(c_l)$.
By the choice of $p$ (in Definition~\ref{definition of R-q}), it follows that all $c_i, c_j, c_k, c_l$ have the same type over $\es$,
and for some $0 < s \leq h$ they all belong to $C_s \setminus C_{s-1}$.

We will prove that $\acl(c_j) = \acl(c_k)$, which implies that
$\acl(c_i) = \acl(c_l)$ and from this we immediately get $R_p(a, c)$.
By symmetry of the argument, it suffices to show that $c_j \in \acl(c_k)$.
By the choice of $c_j$ and $c_k$ we have $p(a, c_j)$ and $p(a, c_k)$ and therefore
$bc_j \equiv bc_k$.
Then Lemma~\ref{when there are independent coordinates with the same type over a}
implies that $c_j \underset{\crd_{s-1}(b)}{\nind} c_k$.
By
Facts~\ref{existence of coordinates}~(iii) and~\ref{rank 1 over its coordinates on the previous level}, 
$\su(c_j / \crd_{s-1}(b)) = 1$ and therefore $c_j \in \acl\big(\{c_k\} \cup \crd_{s-1}(b)\big)$.
By Fact~\ref{fact about C and crd}~(i), $c_j \in \acl(d)$ for some $d \in \{c_k\} \cup \crd_{s-1}(b)$.
As $\su(c_j / \crd_{s-1}(b)) = 1$ we must have $c_j \in \acl(c_k)$.
\hfill $\square$
\\

\begin{lem}\label{acl(c) = acl(a_R)}
Let $a \in M$, $c \in \crd(a)$, $p = \tp(a, c)$ and let $R_p$ be as in Definition~\ref{definition of R-q}.
Then $\acl(c) = \acl(a_{R_p})$.
\end{lem}

\noindent
{\bf Proof.}
By the definition of $R_p$, for every $a' \in [a]_{R_p}$, $c \in \acl(a')$.
Hence 
Lemma~\ref{being in the algebraic closure of an equivalence relation} 
implies that $c \in \acl(a_{R_p})$.
By the definition of $R_p$, $[a]_{R_p}$ is the unique $R_p$-class such that for all $a' \in [a]_{R_p}$,
there is $c'$ with $\acl(c') = \acl(c)$ and $p(a', c')$. Hence $a_{R_p} \in \dcl(c)$.
\hfill $\square$
\\

\noindent
Let $a \in M$. Let $h < \omega$ be minimal such that $a \in \acl(C_h)$.
It follows (from Fact~\ref{fact about C and crd}~(iv)) that $a \in \acl(\crd_h(a))$.

\begin{defin}\label{choice of c-s-i}{\rm
(i) For each $0 < s \leq h$, 
let $\rho_s$ be maximal so that there are $c_{s,1}, \ldots, c_{s,\rho_s} \in \crd_s(a) \setminus \crd_{s-1}(a)$ such that
$\{c_{s,1}, \ldots, c_{s,\rho_s}\}$ is an independent set over $\crd_{s-1}(a)$.
(So $\rho_s$ is the ``dimension'' of $\crd_s(a) \setminus \crd_{s-1}(a)$ over $\crd_{s-1}(a)$.)
We now fix such $c_{s,1}, \ldots, c_{s,\rho_s}$.\\
(ii) For each $0 < s \leq h$ and $1 \leq i \leq \rho_s$, let $p_{s,i} = \tp(a, c_{s,i})$.\\
(iii) For each $0 < s \leq h$ and $1 \leq i \leq \rho_s$, let 
$R_{s, i} = R_{p_{s,i}}$ where $R_{p_{s,i}}$ is like $R_p$ in Definition~\ref{definition of R-q} with $p = p_{s,i}$.
}\end{defin}

\begin{observation}\label{observation about coordinates}{\rm
From Lemma~\ref{when there are independent coordinates with the same type over a}
it follows that, for every $1 \leq s \leq h$ and all
$1 \leq i < j \leq \rho_s$, $p_{s,i} \neq p_{s, j}$. 
And we clearly have $p_{s,i} \neq p_{s',j}$ if $s \neq s'$.
It follows that if $(s, i) \neq (s', i')$ then $R_{s,i}$ is different from $R_{s',i'}$.
}\end{observation}

\begin{defin}\label{ordering of indices}{\rm
Let $I = \{(s, i) : 1 \leq s \leq h, \ 1 \leq i \leq \rho_s\}$ and let `$\preccurlyeq$' be the lexicographic order on $I$,
in other words, $(s, i) \preccurlyeq (s', i')$ if and only if $s < s'$, or $s = s'$ and $i \leq i'$.
}\end{defin}

\noindent
Note that while the ordering in the first coordinate of $(s, i)$ is natural, since $s$ is the ``height'' of $c_{s,i}$,
the order in the second coordinate is arbitrary, since it is given by the arbitrary enumeration $c_{s,1}, \ldots, c_{s,\rho_s}$
of the same elements.

\begin{defin}\label{definition of R-s-i}{\rm
For every $(s,i) \in I$, let 
\[
Q_{s, i} = \bigcap_{(s',i') \preccurlyeq (s,i)} R_{s',i'}.
\]
}\end{defin}

\noindent
Since intersections/conjuctions of equivalence relations are still equivalence relations it follows from
Lemma~\ref{R-q is an equivalence relation} that $Q_{s, i}$ is a $\es$-definable equivalence relation for each $(s,i)$.

\begin{lem}\label{algebraic closure of coordinates and the corresponding equivalence class}
For every $(s,i) \in I$,
\[
\acl(a_{Q_{s,i}}) = \acl\big(\{ c_{s',i'} : (s',i') \preccurlyeq (s,i)\}\big).
\]
\end{lem}

\noindent
{\bf Proof.} 
Let $(s,i) \in I$.
We have $c_{s',i'} \in \acl(a)$ for all $(s',i') \preccurlyeq (s,i)$.
From the definitions of $Q_{s,i}$ and $R_{s',i'}$ it follows that for every $a' \in a_{Q_{s,i}}$, 
$c_{s',i'} \in \acl(a')$ for all $(s',i') \preccurlyeq (s,i)$.
Lemma~\ref{being in the algebraic closure of an equivalence relation}
now implies that 
$c_{s',i'} \in \acl(a_{Q_{s,i}})$ for all $(s',i') \preccurlyeq (s,i)$.

By Lemma~\ref{acl(c) = acl(a_R)},
for every $(s',i') \preccurlyeq (s,i)$,
$\acl(a_{R_{s',i'}}) = \acl(c_{s',i'})$.
From the definition of $Q_{i,s}$ it follows that, for any $a', b' \in M$, $a' \in [b']_{Q_{i,s}}$ if and only if
$a' \in [b']_{R_{s',i'}}$ for all $(s',i') \preccurlyeq (s,i)$.
Consequently,
\[
a_{Q_{s,i}} \in \acl\big(\{a_{R_{s',i'}} : (s',i') \preccurlyeq (s,i) \}\big) = 
\acl\big(\{c_{s',i'} : (s',i') \preccurlyeq (s,i) \}\big). \qquad \qquad \square
\]

\begin{lem}\label{the rank of an element over its predeccessor}
Suppose that $(s,i) \in I$ is not maximal and that $(s',i')$ is the least element in $I$ which is strictly larger 
(with respect to `$\preccurlyeq$') than $(s,i)$.
Then
$\su(a_{Q_{s,i}} / a_{Q_{s',i'}}) =~1$.
\end{lem}

\noindent
{\bf Proof.}
Let $(s,i), (s',i') \in I$ satisfy the assumptions of the lemma.
By Lemma~\ref{algebraic closure of coordinates and the corresponding equivalence class},
it suffices to show that $\su(\bar{c}^+ /\bar{c}) = 1$, where
\[
\bar{c} = \big(c_{t,j} : (t,j) \preccurlyeq (s,i)\big) \ \text{ and } \
\bar{c}^+ = \big(c_{t,j} : (t,j) \preccurlyeq (s',i')\big).
\]
To show this we only need to show that  $\su(c_{s',i'} / \bar{c}) = 1$.

We have two cases.
First, suppose that $s = s'$. Then $i' = i+1$.
By the choice of the elements $c_{t,j}$ and
Facts~\ref{existence of coordinates}~(iii) and~\ref{rank 1 over its coordinates on the previous level}, 
we get
$\su\big(c_{s,i+1}/ \{c_{t,j} : (t,j) \preccurlyeq (s-1,\rho_{s-1})\}\big) =
\su\big(c_{s,i+1}/ \crd_{s-1}(a)\big) = 1$.
And we also have that $\{c_{s,1}, \ldots, c_{s,i+1}\}$ is independent over $\crd_{s-1}(a)$.
Therefore, $\su(c_{s,i+1} / \bar{c}) = 1$.

Now suppose that $s' = s+1$, so $i = \rho_s$ and $i' = 1$.
As in the previous case we get
$\su\big(c_{s+1,1} / \bar{c}) =   \su\big(c_{s+1,1} / \{c_{t,j} : (t,j) \preccurlyeq (s,\rho_{s})\}\big) 
= \su\big(c_{s+1,1}/ \crd_s(a)\big) = 1$ and we are done.
\hfill $\square$
\\

\noindent
Now we can finish the proof of part~(i) of Theorem~\ref{theorem coordinatizing M by equivalence classes}.
Recall that (by Fact~\ref{fact about C and crd}~(iv)) $\acl(a) = \acl(\crd_h(a))$ and therefore 
(using Lemma~\ref{algebraic closure of coordinates and the corresponding equivalence class})
\[
a \in \acl\big(\{c_{s,i} : (s,i) \in I\}) = \acl\big(\{a_{Q_{s,i}} : (s,i) \in I\}\big).
\]
Since $Q_{s,i}$ refines $Q_{s',i'}$ if $(s',i') \preccurlyeq (s,i)$ we get
$a \in \acl(a_{Q_{h,\rho_h}})$.
Since $c_{1,1} \in C_1$ we have 
(using Lemma~\ref{algebraic closure of coordinates and the corresponding equivalence class}
and Fact~\ref{existence of coordinates}~(iii))
$\su(a_{Q_{1,1}}) = \su(c_{1,1}) = 1$.
From this and Lemma~\ref{the rank of an element over its predeccessor} it follows, via the Lascar equation,
that $\su(a) = |I|$.
Thus the sequence of $\es$-definable equivalence relations that we are looking for is,
again using Lemma~\ref{the rank of an element over its predeccessor},
\[
(Q_{s,i} : (s,i) \in I),
\]
ordered by `$\preccurlyeq$'.

\subsection{Proof of parts~(ii) and~(iii) of Theorem~\ref{theorem coordinatizing M by equivalence classes}}
\label{proof of parts ii and iii}

The assumptions and framework in this subsection are the same as in  the previous
(i.e. the proof of part~(i)).

Suppose that $a, b \in M$, $\bar{c} \in C$ and $a \underset{\bar{c}}{\nind} b$.
By Lemma~\ref{a more general independence criterium},
there is $d \in \crd(a) \cap \crd(b)$ such that $d \notin \acl(\bar{c})$.
Let $p = \tp(a, d)$ and let $R = R_p$ be as in
Definition~\ref{definition of R-q}.
By Lemma~\ref{R-q is an equivalence relation},
$R$ is a $\es$-definable equivalence relation.
By Lemma~\ref{acl(c) = acl(a_R)},
$\acl(d) = \acl(a_{R})$.
Since $d \in \crd(b) \subseteq \acl(b)$ we get $a_{R} \in \acl(b)$.
By assumption, $d \notin \acl(\bar{c})$ and hence $a_{R} \notin \acl(\bar{c})$.
Then there are distinct $a'_i \in M\meq$, for $i < \omega$, such that
$\tp(a'_i / \bar{c}) = \tp(a_R / \bar{c})$, for all $i < \omega$.
The type $\tp(a / a_R)$ contains a formula, $\varphi(x, a_R)$ which expresses that
\begin{itemize}
\item[] ``$x$ belongs to the equivalence class (represented by) $a_R$''.
\end{itemize}
Since $\{\varphi(x, a'_i) : i < \omega\}$ is clearly 2-inconsistent it follows that 
$a \underset{\bar{c}}{\nind} a_R$.
This concludes the proof of Theorem~\ref{theorem coordinatizing M by equivalence classes}~(ii).

Now assume, in addition, that every binary $\es$-definable relation on $M$ is symmetric.
Suppose that $a, b, \bar{c} \in M$ and $a \underset{\bar{c}}{\nind} b$.
Just as in the proof of part~(ii), we get $d \in \crd(a) \cap \crd(b)$ such that $d \notin \acl(\bar{c})$.
By letting $p = \tp(a, d)$ and $R = R_p$ be just as in the proof of part~(ii), we conclude (just as in part~(ii))
that 
$R$ is a $\es$-definable equivalence relation and $\acl(d) = \acl(a_{R})$.

Let $0 < s \leq r$ be such that $d \in C_s \setminus C_{s-1}$.
By Fact~\ref{fact about C and crd}, $C_s \setminus C_{s-1}$ is a trivial pregeometry.
So if there is $e \in C_s \setminus C_{s-1}$ such that $\acl(d) = \acl(e)$
and $p(b, e)$, then $R_p(a, b)$ (by Definition~\ref{definition of R-q}) so $R(a, b)$ and hence $a_R = b_R$.
Since $\acl(d) = \acl(a_R)$ and $d \notin \acl(\bar{c})$ we must have $a_R \notin \acl(\bar{c})$.

Now suppose (towards a contradiciton) that, for every $e \in  C_s \setminus C_{s-1}$ such that $\acl(d) = \acl(e)$,
we have $\tp(b, e) \neq p$.

Let $a'$ realize a nondividing extension of 
\[
\tp\big(a \big/ (\crd_{s-1}(a) \cap \crd_{s-1}(b))  \cup \{d\}\big) \ \text{ to } \
( \crd_{s-1}(a) \cap \crd_{s-1}(b)) \cup \{d, b\}.
\]
Then $a'$ is independent from $b$ over $\big( \crd_{s-1}(a') \cap \crd_{s-1}(b)\big) \cup \{d\}$.
As $C_s \setminus C_{s-1}$ is a trivial pregeometry, it follows that 
if $e \in \big(\crd_s(a') \cap \crd_s(b)\big) \setminus C_{s-1}$, then $\acl(e) = \acl(d)$.
By assumption, 
for every $e \in \big(\crd_s(a') \cap \crd_s(b)\big) \setminus C_{s-1}$, $\tp(b, e) \neq p$.
But then $\tp(a', b) \neq \tp(b, a')$. Since $a', b \in M$ we get $\tp_\mcM(a', b) \neq \tp_\mcM(b, a')$.
As every complete type over $\es$ is isolated it follows that there is a binary $\es$-definable relation 
which is not symmetric, which contradicts an assumption of part~(iii). Thus the proof of part~(iii) is finished.

\subsection{Proof of Corollaries~\ref{information about nonbinary structures}
and~\ref{information about nonbinary primitive structures}}
\label{proofs of two corollaries}

Suppose that $\mcM$ is $\omega$-categorical, supersimple with finite SU-rank and with trivial dependence.
Moreover, suppose that part~(i) of Theorem~\ref{theorem coordinatizing M by equivalence classes}
does not hold for $\mcM$.
The proof of part~(i) of Theorem~\ref{theorem coordinatizing M by equivalence classes} only uses
\begin{itemize}
\item results from Section~\ref{Coordinatization} 
all of which hold for all $\omega$-categorical, supersimple structures with finite SU-rank and with trivial dependence,
\item Lemma~\ref{when there are independent coordinates with the same type over a}, and
\item results from Section~\ref{proof of part i} which,
besides Lemma~\ref{when there are independent coordinates with the same type over a}
only depend on the assumption that $\mcM$ is 
$\omega$-categorical, supersimple with finite SU-rank and with trivial dependence.
\end{itemize}
So, assuming that part~(i) of Theorem~\ref{theorem coordinatizing M by equivalence classes} fails for $\mcM$, it must be
because Lemma~\ref{when there are independent coordinates with the same type over a} fails for $\mcM$.
But the proof of Lemma~\ref{when there are independent coordinates with the same type over a}
is a proof by contradiction. It assumes that 
Lemma~\ref{different coordinates have different types over acl(emptyset)} holds 
(and consequently Lemma~\ref{about coordinates and equivalence classes} holds) and that
Lemma~\ref{when there are independent coordinates with the same type over a} fails, 
and then finds $a^*, b^*, a', b', b'' \in M$ such that
\[
a^*b^*b' \equiv a^*b^*b'' \ \text{ and } \ a'b' \equiv a'b'', \ \text{ but } \
a^*b^*a'b' \not\equiv a^*b^*a'b''.
\]
This finishes the proof of Theorem~\ref{information about nonbinary structures}
unless Lemma~\ref{different coordinates have different types over acl(emptyset)} fails for $\mcM$.
But if  Lemma~\ref{different coordinates have different types over acl(emptyset)} fails,
then (by its proof) there are $a, a', a'', a^* \in M$ such that
$aa'' \equiv aa^*$ and $a'a'' \equiv a'a^*$, but $aa'a'' \not\equiv aa'a^*$.
This finishes the proof of Theorem~\ref{information about nonbinary structures}.

Now we prove Corollary~\ref{information about nonbinary primitive structures}.
Suppose that $\mcM$ is $\omega$-categorical, supersimple with finite SU-rank and with trivial dependence.
Moreover, assume that $\mcM$ has no $\es$-definable equivalence relation on $M$
with infinitely many infinite equivalence classes.
By the proof of \cite[Lemma~3.3]{Kop-one-based}, $M \subseteq \acl(C_1)$.\footnote{
Lemma~3.3 in \cite{Kop-one-based} assumes that $\mcM$ is primitive, but its proof only needs the assumption that
that there is no $\es$-definable equivalence relation on $M$ which  has infinitely many infinite equivalence classes.}
Furthermore, assume that $\su(\mcM) > 1$, so $\su(a) > 1$ for some $a \in M$.

Suppose that $\su(a) = \rho_1 > 1$.
By Fact~\ref{fact about C and crd}~(iv),
$\acl(a) = \acl(\crd_1(a))$, and hence there are $c_{1,1}, \ldots, c_{1,\rho_1} \in \crd_1(a)$
such that $\{c_{1,1}, \ldots, c_{1,\rho_1}\}$ is an independent set over $\es$ and 
$\acl(a) = \acl(c_{1,1}, \ldots, c_{1, \rho_1})$.

Suppose that Lemma~\ref{when there are independent coordinates with the same type over a} holds for $\mcM$.
Then
$ac_{1,i} \not\equiv ac_{1,j}$ whenever $i \neq j$.
Let $p = \tp(a, c_{1,1})$.
By Lemma~\ref{R-q is an equivalence relation},
$R_p$ (as in Definition~\ref{definition of R-q}) 
is a $\es$-definable equivalence relation on $M$.
Since $\su(c_{1,1}) = 1$, it follows from
Lemma~\ref{the rank of an element over its predeccessor}
that $R_p$ has infinitely many infinite equivalence classes.
This contradicts the assumptions of
Corollary~\ref{information about nonbinary primitive structures}.

Hence Lemma~\ref{when there are independent coordinates with the same type over a} fails for $\mcM$.
Then, in the same way as in the proof of Corollary~\ref{information about nonbinary structures},
we find $a_1, a_2, a_3, a_4, b_1, b_2, b_3, b_4 \in M$ such that
$a_ia_j \equiv b_ib_j$ for all $i$ and $j$, but
$a_1a_2a_3a_4 \not\equiv b_1b_2b_3b_4$.
This completes the proof of Corollary~\ref{information about nonbinary primitive structures}.

\section{Extension properties}\label{Extension properties}

\noindent
We are interested in knowing under what conditions two or more types are subtypes of a single type.
More precisely, if $\bar{a}_i, \bar{b}_i \in M$, for $i = 1, \ldots, n$, under what circumstances is there $\bar{a} \in M$
such that $\tp(\bar{a}, \bar{b}_i) = \tp(\bar{a}_i, \bar{b}_i)$ for all $i = 1, \ldots, n$?
Under rather general conditions, the
answer is yes for the Rado graph, 
the ``generic bipartite graph'', as well as a number of other structures that can
be constructed by procedures that involve a ``high degree of randomness''.\footnote{
The most up to date study of extension problems in the context of binary $\omega$-categorical structures
is probably  \cite{Ahl}, by Ahlman, where more references can be found.}
Therefore, the idea here is that if the answer is `yes' under fairly general conditions,
then this is a manifestation of a ``high degree of randomness''.

\begin{defin}\label{definition of extension problem}{\rm
Here we call the following an {\em extension problem} of $\mcM$:
\begin{itemize}
\item[] Suppose that $\bar{a}_1, \ldots, \bar{a}_n, \bar{b}_1, \ldots, \bar{b}_n \in M$. 
We ask: is there $\bar{e} \in M$ such that
$\tp(\bar{e}, \bar{b}_i) = \tp(\bar{a}_i, \bar{b}_i)$ for all $i = 1, \ldots, n$?
If such $\bar{e} \in M$ exists then we say that {\em the extension problem of 
$\tp(\bar{a}_i, \bar{b}_i)$, $i = 1, \ldots, n$, has a solution} and call $\bar{e}$ a {\em solution} to this extension problem.
\end{itemize}
}\end{defin}

\noindent
Observe that since we will assume that $\mcM$ is homogeneous (hence $\omega$-saturated)
it follows that if an extension problem has a solution $\bar{e}$ in some elementary extension of $\mcM$,
then it also has a solution in $\mcM$.

Note also that if we have $\bar{a}_i$ and $\bar{b}_i$ as above and, for every $i = 1, \ldots, n$, there is
$\bar{a}'_i$ such that for every $i < n$, 
$\tp(\bar{a}'_{i+1}, \bar{b}_1, \ldots, \bar{b}_i) = \tp(\bar{a}'_i, \bar{b}_1, \ldots, \bar{b}_i)$
and $\tp(\bar{a}'_{i+1}, \bar{b}_{i+1}) = \tp(\bar{a}_{i+1}, \bar{b}_{i+1})$, then we
have $\tp(\bar{a}'_n, \bar{b}_i) = \tp(\bar{a}_i, \bar{b}_i)$ for every $i = 1, \ldots, n$.
Therefore we will only consider the problem of extending two types.

Here we study binary relational structures with elimination of quantifiers.
Under this assumption, if $\bar{c} = (c_1, \ldots, c_k)$, and $\bar{e}$ is a solution to
the extension problem of the types
$\tp(\bar{a}, c_1), \ldots, \tp(\bar{a}, c_k), \tp(\bar{b}, \bar{d})$,
then $\bar{e}$ is also a solution to the extension problem of the types
$\tp(\bar{a}, \bar{c})$ and $\tp(\bar{b}, \bar{d})$.
And as pointed out above, the extension problem of the types
$\tp(\bar{a}, c_1), \ldots, \tp(\bar{a}, c_k), \tp(\bar{b}, \bar{d})$ can be reduced to a sequence of $k$
extension problems of two types of the form $\tp(\bar{a}', c')$ and $\tp(\bar{b}', \bar{d}')$, where $c'$ is a single element.

By considering one coordinate at a time in the sequences $\bar{a}_1, \ldots, \bar{a}_n$, and using our observations above,
it follows that
the extension problem of
the types $\tp(\bar{a}_1, \bar{b}_1), \ldots, \tp(\bar{a}_n, \bar{b}_n)$
can be reduced to a sequence of extension problems of two types of the form
$\tp(a', c')$ and $\tp(b', \bar{d}')$, where $a', b'$ and $c'$ are single elements.
Therefore we will only consider the extension problem of
two types $\tp(a, c)$ and $\tp(b, \bar{d})$, where $a, b$ and $c$ are a single elements.
Recall that notation~\ref{remark about acl and dcl} is in effect in this section.

\begin{theor}\label{main result regarding near randomness}
Suppose that $\mcM$ is binary, simple and homogeneous.
Let $a, b, c, \bar{d} \in M$.
\begin{itemize}
\item[(i)] There is a $\es$-definable equivalence relation $R$ on $M$ such that
$c \underset{c_R}{\ind} \bar{d}$.
\item[(ii)] If for some $R$ as in part~(i),
\begin{itemize}
\item[] $a \underset{c_R}{\ind} c$, \ $b \underset{c_R}{\ind} \bar{d}$ \ and \
$\tp(a / \acl(c_R)) = \tp(b / \acl(c_R))$,
\end{itemize}
then the extension problem of $\tp(a, c)$ and $\tp(b, \bar{d})$ has a solution. Otherwise it may not have a solution,
not even when $\bar{d}$ is a single element.
\end{itemize}
\end{theor}

\noindent
{\bf Proof.}
This follows from Lemmas~\ref{finding an equivalence class over which c and d are independent}~--~\ref{when there are no obstacles}
(and the examples in Sections~\ref{Cross cutting equivalence relations}~--~\ref{omega-pedes}).
\hfill $\square$
\\

\noindent
{\bf \em In the rest of this section we assume that $\mcM$ is binary, simple and homogeneous,
so Theorem~\ref{theorem coordinatizing M by equivalence classes} applies.}

\begin{lem}\label{finding an equivalence class over which c and d are independent}
For all $c, \bar{d} \in M$ there is a $\es$-definable equivalence relation $R$ such that $c \underset{c_R}{\ind} \bar{d}$.
\end{lem}

\noindent
{\bf Proof.}
Recall that, by Assumption~\ref{assumption about C},
$M \subseteq U \subseteq M\meq$ and only finitely many sorts are represented in $U$.
The only assumption on $U$ that is necessary for Fact~\ref{existence of coordinates} to hold
is that only finitely many sorts are represented in $U$ .
Since there are only finitely many $\es$-definable equivalence relations on $M$, we may, without loss of generality, 
assume that for every $\es$-definable equivalence relation $E$ on $M$ and every $a \in M$, $a_E \in U$ and hence $a_E \in C$.

Now we prove~(i).
Let $c, \bar{d} \in M$.
If $c \ind \bar{d}$ then we can take $R$ to be the equivalence relation with only one equivalence class.
So suppose that $c \nind \bar{d}$. 
Then $c \nind d$ for some $d \in \bar{d}$. 
By Theorem~\ref{theorem coordinatizing M by equivalence classes}~(ii),
there is a $\es$-definable equivalence relation $R_1$ such that $c \nind c_{R_1}$ and $c_{R_1} \in \acl(d) \subseteq \acl(\bar{d})$.
If $c \underset{c_{R_1}}{\ind} \bar{d}$ then we are done with $R = R_1$.
If not, then $c \underset{c_{R_1}}{\nind} d$ for some $d \in \bar{d}$, and
(by Theorem~\ref{theorem coordinatizing M by equivalence classes}~(ii)\footnote{
Here we use that $c_{R_1} \in C$, which is why we need the argument in the first paragraph of the proof.}),
there is a $\es$-definable equivalence relation $R_2$ such that $c \underset{c_{R_1}}{\nind} c_{R_2}$
and $c_{R_2} \in \acl(d) \subseteq \acl(\bar{d})$.
If $c \underset{c_{R_2}}{\ind} \bar{d}$ then we are done with $R = R_2$.
If not, we continue in the same way.
Since $\mcM$ has finite SU-rank we will, after finitely many iterations of this procedure,
find a $\es$-definable equivalence relation $R_k$ such that $a \underset{c_{R_k}}{\ind} \bar{d}$. 
(Or alternatively, one could appeal to the fact that there are only finitely many $\es$-definable equivalence relations on $M$.)
\hfill $\square$

\begin{lem}\label{the first obstacle}
Suppose that $a, b, c, \bar{d} \in M$ and that $R$ is a $\es$-definable equivalence relation on $M$ such that
$c \underset{c_R}{\ind} \bar{d}$.
If $a \underset{c_R}{\nind} c$ or $b \underset{c_R}{\nind} \bar{d}$, 
then the extension problem of $\tp(a, c)$ and $\tp(b, \bar{d})$ may {\em not} have a solution.
\end{lem}

\noindent
{\bf Proof.}
This is shown by the examples in Sections~\ref{Cross cutting equivalence relations}
and~\ref{Bipedes with bicoloured legs}.
\hfill $\square$

\begin{lem}\label{the second obstacle}
Suppose that $a, b, c, \bar{d} \in M$, that $R$ is a $\es$-definable equivalence relation on $M$ such that
$c \underset{c_R}{\ind} \bar{d}$ and that $a \underset{c_R}{\ind} c$ and $b \underset{c_R}{\ind} d$.
If $\tp(a / \acl(c_R)) \neq \tp(b / \acl(c_R))$ then the extension problem 
of $\tp(a, c)$ and $\tp(b, \bar{d})$ may not have a solution.
\end{lem}

\noindent
{\bf Proof.}
This is shown by the example in Section~\ref{omega-pedes}.
\hfill $\square$

\begin{lem}\label{when there are no obstacles}
Suppose that $a, b, c, \bar{d} \in M$, $R$ is a $\es$-definable equivalence relation on $M$ such that
$c \underset{c_R}{\ind} \bar{d}$, $a \underset{c_R}{\ind} c$, $b \underset{c_R}{\ind} d$ and
$\tp(a / \acl(c_R)) = \tp(b / \acl(c_R))$.
Then the extension problem of $\tp(a, c)$ and $\tp(b, \bar{d})$ has a solution.
\end{lem}

\noindent
{\bf Proof.}
If the premisses of the lemma are satisfied, then all premisses of the independence theorem of simple theories are satisfied, 
and hence a solution exists in some elementary extension of $\mcM$. Since $\mcM$ is $\omega$-saturated we
find a solution in $M$.
\hfill $\square$

\section{Examples}

\noindent
In sections~\ref{Cross cutting equivalence relations}~--~\ref{omega-pedes}
we give examples that prove the claims made in 
Remark~\ref{remarks on the coordinatization by equivalence relations}
and in
lemmas~\ref{the first obstacle}
and~\ref{the second obstacle}.
Section~\ref{Metric spaces} tells how certain metric spaces
fit nicely into the context of this article when viewed as binary structures
(namely, $\mcR$-Uryshon spaces for finite distance monoids $\mcR$).

\subsection{Cross cutting equivalence relations}\label{Cross cutting equivalence relations}

In this subsection we prove Lemma~\ref{the first obstacle}.
This is also done, in a stronger sense, in Section~\ref{Bipedes with bicoloured legs},
but the example of this section may nevertheless be instructive because of its simplicity.

Let $\mcM = (M, P^\mcM, Q^\mcM)$, where $M$ is a countably infinite set and
$P^\mcM$ and $Q^\mcM$ are equivalence relations such that the equivalence relation $P^\mcM \cap Q^\mcM$
\begin{itemize}
\item partitions every equivalence class of $P^\mcM$ into infinitely many 
parts, all of which are infinite, and
\item partitions every equivalence class of $Q^\mcM$ into infinitely many 
parts, all of which are infinite.
\end{itemize}
It is a basic exercise to show that $\mcM$ is homogeneous and stable 
(hence $\omega$-stable\footnote{
By the work on stable homogeneous structures. See for instance~\cite{Lach97}.} 
and consequently superstable) with SU-rank 2.

Let $X_1$ and $X_2$ be two distinct equivalence classes of $P^\mcM$ and let
$Y_1$ and $Y_2$ be two distinct equivalence classes of $Q^\mcM$.
Pick $a \in X_1 \cap Y_1$, $b \in X_2 \cap Y_1$, $c \in X_1 \cap Y_2$ and $d \in  X_2 \cap Y_2$.
Then it is straightforward to verify that $c \underset{c_Q}{\ind} d$, 
$a \underset{c_Q}{\nind} c$ and $b \underset{c_Q}{\nind} d$, where
`$c_Q$' is shorthand for `$c_{Q^\mcM}$'.
Moreover, the extension problem of $\tp_\mcM(a, c)$ and $\tp_\mcM(b, d)$ does not have a solution,
because if $e$ would be a solution then $\mcM \models P(e, c) \wedge P(e, d)$, so
$\mcM \models P(c, d)$, contradicting the choice of $c$ and $d$.

\subsection{Bipedes with bicoloured legs}\label{Bipedes with bicoloured legs}

In this subsection we prove the claims made in 
Remark~\ref{remarks on the coordinatization by equivalence relations} and Lemma~\ref{the first obstacle}.
For any set $A$, let $[A]^2 = \{X \subseteq A : |X| = 2\}$.
Let 
\[
\mcN^- = (\mbbN \cup [\mbbN]^2, F^{\mcN^-}, L^{\mcN^-}),
\]
where 
\[
F^{\mcN^-} = \mbbN \ \text{ and } \ 
L^{\mcN^-} = \big\{ (\{m, n\}, k) : \{m, n\} \in [N]^2 \text{ and } k \in \{m, n\} \big\}.
\]
We can think of the elements of $F^{\mcN^-} = \mbbN$ as ``feet'' and elements of $[\mbbN]^2$ as ``bodies''.
Each body $\{m, n\} \in [\mbbN]^2$ has two feet, namely $m$ and $n$.
Clearly, some different bodies, like $\{1, 2\}$ and $\{2, 3\}$, share a foot, while others do not.
We can also imagine any given pair $(\{m, n\}, n) \in L^{\mcN^-}$ as a ``leg'' which joins the
body $\{m, n\}$ to the foot $n$.
We further imagine that for every body, one of its legs is coloured ``blue'' and the other is coloured ``red''.
Moreover, the decision regarding which one is blue and which one is red is taken randomly and independently of
the colouring of the legs of other ``bodies''. Note that only legs are coloured. A given foot may be the end of a blue leg and also the end
of red leg, in which case the later leg belongs to another body than the first leg.

More formally, we construct such a structure as follows.
Let $B$ and $R$ (for ``blue'' and ``red'') be new binary relation symbols and 
let $\Omega$ be the set of expansions 
\[
\mcN = (\mbbN \cup [\mbbN]^2, F^\mcN, L^\mcN, B^\mcN, R^\mcN) 
\]
of $\mcN^-$ which satisfy the following sentences:
\begin{align*}
&\forall x, y \Big( \big[ B(x,y) \vee R(x,y) \big] \rightarrow L(x,y) \Big), \\
&\forall x, y \Big( L(x,y) \rightarrow \big[ (B(x,y) \wedge \neg R(x,y)) \vee (R(x,y) \wedge \neg B(x,y)) \big] \Big),
\ \text{ and } \\
&\forall x \Big( \neg F(x) \rightarrow \exists y, z \big[ B(x,y) \wedge R(x,z) \big] \Big).
\end{align*}

\noindent
For any set $X$ let $2^X$ denote the set of functions from $X$ to $\{0,1\}$
For every finite $A \subseteq \mbbN$ and every $f \in 2^A$, 
let $\langle A, f \rangle = \{g \in 2^\mbbN : g(n) = f(n) \text{ for all } n \in A\}$.
If $|A| = m$ then we let $\mu_0(\langle A, f \rangle) = 2^{-m}$.
By standard notions and results in measure theory, there is a $\sigma$-algebra $\Sigma \subseteq 2^\mbbN$,
containing all $\langle A, f \rangle$ for finite $A$ and $f \in 2^A$, and a 
countably subadditive probability measure $\mu : \Sigma \to \mbbR$ which
extends $\mu_0$.\footnote{These notions and results can be found in, for example, \cite[Chapters~1.1--1.4]{Fri}.}
Let $\lambda : [\mbbN]^2 \to \mbbN$ be a bijection.
For every $f \in 2^\mbbN$ we get an expansion $\mcN_f \in \Omega$ of $\mcN$ that satisfies:
\begin{itemize}
\item[] For every $\{m, n\} \in [\mbbN]^2$ with $m < n$,
if $f(\lambda(\{m, n\})) = 0$, then $\mcN_f \models B(\{m, n\}, m) \wedge R(\{m, n\}, n)$, 
and otherwise $\mcN_f \models R(\{m, n\}, m) \wedge B(\{m, n\}, n)$.
\end{itemize}
Moreover, it is clear that for every $\mcN \in \Omega$ there is a unique $f \in 2^\mbbN$ such that $\mcN = \mcN_f$.
Via this bijection between $2^\mbbN$ and $\Omega$ we may also view $\Omega$ as a probability space.

\begin{lem}\label{N has extension properties}
There is $\mcN \in \Omega$ with the following property.
Let $0 < n < \omega$, $a_1, \ldots, a_n \in \mbbN$ and $f : \{1, \ldots, n\} \to \{0,1\}$.
Then there are distinct $b_i \in \mbbN \setminus \{a_1, \ldots, a_n\}$, for all $i < \omega$, such that, 
for every $i < \omega$ and every $1 \leq k \leq n$, the following holds:
\begin{itemize}
\item If $f(k) = 0$ then $\mcN \models B(\{a_k, b_i\}, b_i) \wedge R(\{a_k, b_i\}, a_k)$.
\item If $f(k) = 1$  then $\mcN \models R(\{a_k, b_i\}, b_i) \wedge B(\{a_k, b_i\}, a_k)$.
\end{itemize}
\end{lem}

\noindent
{\bf Proof.}
We will prove that with probability 1 a structure in $\Omega$ has the stated property.
By countable subadditivity of $\mu$, it suffices to show the following:
\begin{itemize}
\item[] For any choice of $0 < n < \omega$, $a_1, \ldots, a_n \in \mbbN$, $f : \{1, \ldots, n\} \to \{0,1\}$
and distinct $b_j^i \in \mbbN \setminus \{a_1, \ldots, a_n\}$ for $i, j < \omega$,
\begin{align*}
&\mu(X_i) = 0, \text{ for every $i < \omega$, where}  \\
&X_i = \Big\{g \in 2^\mbbN : \text{ for all $j < \omega$ there is $1 \leq k \leq n$ such that } \\
&g(\lambda(\{a_k, b_j^i\})) \neq f(\lambda(\{a_k, b_j^i\})) \Big\}.
\end{align*}
\end{itemize}
By using the definition of $\mu_0$ and the fact that $\mu$ extends $\mu_0$ we get
\[
\mu(X_i) \ \leq \ \big(1 - 2^{-(n+1)}\big)^j 
\]
for every $i < \omega$ and every $j < \omega$. Hence $\mu(X_i) = 0$ for every $i < \omega$ and the proof is finished.
\hfill $\square$
\\

\noindent
{\bf \em For the rest of this subsection we assume that $\mcN$ is like in Lemma~\ref{N has extension properties}}.

\begin{defin}\label{definition of closure in N}{\rm
(i) For every $A \subseteq \mbbN \cup [\mbbN]^2$, 
$\cl'(A) = A \cup \{b \in \mbbN : \exists a \in A \cap [\mbbN]^2, b \in a\}$.\\
(ii) For every $A \subseteq \mbbN \cup [\mbbN]^2$, 
$\cl''(A) = A \cup \{b \in [\mbbN]^2 : \exists m, n \in A \cap \mbbN, b = \{m, n\} \}$.\\
(iii) For every $A \subseteq \mbbN \cup [\mbbN]^2$, 
$\cl(A) = \cl''(\cl'(A))$.\\
(iv) We say that $A \subseteq \mbbN \cup [\mbbN]^2$ is {\em closed} if $\cl(A) = A$.
}\end{defin}

\begin{lem}\label{cl(A) is a subset of dcl(A)}
Suppose that $A \subseteq \mbbN \cup [\mbbN]^2$ and $a \in \cl(A)$.
Then there is $B \subseteq A$ such that $|B| \leq 2$ and $a \in \dcl_\mcN(B)$.
If $a \in \mbbN$, then there is $b \in A$ such that $a \in \dcl_\mcN(b)$.
\end{lem}

\noindent
{\bf Proof.}
This is because,
\begin{itemize}
\item[(a)] for any two (different) feet there is a unique body which has precisely these two feet, and
\item[(b)] every body has a unique foot on the other end of its blue leg and a unique foot on the other end of its red leg.
\hfill $\square$
\end{itemize}

\begin{lem}\label{extensions of closed sets by elements in F}
Suppose that $\{a_1, \ldots, a_n\}, \{ b_1, \ldots, b_n\} \in \mbbN \cup [\mbbN]^2$ are two
closed sets such that $(a_1, \ldots, a_n) \equiv^{at}_\mcN (b_1, \ldots, b_n)$.\\
(i) For every $a_{n+1} \in \mbbN \cup [\mbbN]^2$ there is $b_{n+1} \in \mbbN \cup [\mbbN]^2$ such that
$\cl(a_1, \ldots, a_{n+1}) \setminus \{a_1, \ldots, a_{n+1}\}$ and
$\cl(b_1, \ldots, b_{n+1}) \setminus \{b_1, \ldots, b_{n+1}\}$
can be enumerated as $a'_1, \ldots, a'_m$ and $b'_1, \ldots, b'_m$, respectively,
so that
\[
(a_1, \ldots, a_{n+1}, a'_1, \ldots, a'_m) \equiv^{at}_\mcN
(b_1, \ldots, b_{n+1}, b'_1, \ldots, b'_m).
\]
(ii) There is an automorphism $\sigma$ of $\mcN$ such that $\sigma(a_i) = b_i$ for every $1 \leq i \leq n$.\\
(iii) $\mcN$ is $\omega$-categorical.\\
(iv) If $a, b \in \mbbN$ or if $a, b \in [\mbbN]^2$, then $\tp_\mcN(a) = \tp_\mcN(b)$.
\end{lem}

\noindent
{\bf Proof.}
(i) We consider two cases. First assume that $a_{n+1} \in \mbbN \setminus \{a_1, \ldots, a_n\}$.
Without loss of generality, assume that
$\{a_1, \ldots, a_n\} \cap \mbbN = \{a_1, \ldots, a_k\}$ for some $k \leq n$.
Then $\{b_1, \ldots, b_n\} \cap \mbbN = \{b_1, \ldots, b_k\}$
Since $(a_1, \ldots, a_n) \equiv^{at}_\mcN (b_1, \ldots, b_n)$ it suffices to find $b_{n+1} \in \mbbN$
such that for every $1 \leq i \leq k$:
\begin{itemize}
\item If $B(\{a_{n+1}, a_i\}, a_{n+1})$ then $B(\{b_{n+1}, b_i), b_{n+1})$.
\item If $R(\{a_{n+1}, a_i\}, a_{n+1})$ then $R(\{b_{n+1}, b_i), b_{n+1})$.
\end{itemize}
But Lemma~\ref{N has extension properties} guarantees that such $b_{n+1} \in \mbbN$ exists.

Now suppose that $a_{n+1} = \{i, j\} \in [\mbbN]^2  \setminus \{a_1, \ldots, a_n\}$.
Then at least one of $i$ or $j$ does not belong to $\{a_1, \ldots, a_n\}$,
because this set is, by assumption, closed.
First, suppose that $i \in \{a_1, \ldots, a_n\}$ and $j \notin \{a_1, \ldots, a_n\}$.
Without loss of generality, assume that $i = a_1$.
Then, by the previous case, we find $j' \in \mbbN$ such that
$\cl(a_1, \ldots, a_n, j) \setminus \{a_1, \ldots, a_n, j\}$ and
$\cl(b_1, \ldots, b_n, j') \setminus \{b_1, \ldots, b_n, j'\}$
can be enumerated as $a'_1, \ldots, a'_m$ and $b'_1, \ldots, b'_m$, respectively,
so that
\[
(a_1, \ldots, a_n, j, a'_1, \ldots, a'_m) \equiv^{at}_\mcN
(b_1, \ldots, b_n, j', b'_1, \ldots, b'_m).
\]
Moreover, since these sequences are closed, there is $1 \leq l \leq m$ such that
$a_{n+1} = \{i,j\} = \{a_1, j\} = a'_l$ and hence $\{b_1, j'\} = b'_l$, so we are done
by taking $b_{n+1} = b'_l$.

Now suppose that $i, j \notin \{a_1, \ldots, a_n\}$.
Then we apply what we have already proved twice.
First we find we find $i' \in \mbbN$ such that 
$\cl(a_1, \ldots, a_n, i) \setminus \{a_1, \ldots, a_n, i\}$ and
$\cl(b_1, \ldots, b_n, i') \setminus \{b_1, \ldots, b_n, i'\}$
can be enumerated as $a'_1, \ldots, a'_m$ and $b'_1, \ldots, b'_m$, respectively,
so that
\[
(a_1, \ldots, a_n, i, a'_1, \ldots, a'_m) \equiv^{at}_\mcN
(b_1, \ldots, b_n, i', b'_1, \ldots, b'_m).
\]
Then we find $j' \in \mbbN$ such that
$\cl(a_1, \ldots, a_n, i, a'_1, \ldots, a'_m, j) \setminus \{a_1, \ldots, a_n, i, a'_1, \ldots, a'_m, j\}$ and
$\cl(b_1, \ldots, b_n, i', b'_1, \ldots, b'_m, j') \setminus \{b_1, \ldots, b_n, i',  b'_1, \ldots, b'_m, j'\}$
can be enumerated as $a''_1, \ldots, a''_s$ and $b''_1, \ldots, b''_s$, respectively,
so that
\[
(a_1, \ldots, a_n, i, a'_1, \ldots, a'_m, j, a''_1, \ldots, a''_s) \equiv^{at}_\mcN
(b_1, \ldots, b_n, i', b'_1, \ldots, b'_m, j', b''_1, \ldots, b''_s).
\]
Then $a_{n+1} = \{i, j\} = a''_l$ for some $l$, and we take $b_{n+1} = \{i', j'\} = b''_l$.

(ii) By part~(i), we can carry out a standard back and forth argument to produce an automorphism $f$ such
that $f(a_i) = b_i$ for all $i$.

(iii) By the definition of `$\cl$' it is clear that, for every finite $A \subseteq \mbbN \cup [\mbbN]^2$, 
$|\cl(A)| \leq 3|A| + \binom{3|A|}{2}$.
Together with part~(ii) this implies that there are, up to equivalence in $Th(\mcN)$, only finitely many formulas
with free variables $x_1, \ldots, x_n$, for every $n < \omega$. Hence $\mcN$ is $\omega$-categorical.

(iv) If $a, b \in \mbbN$, then $\{a\}$ and $\{b\}$ are closed and $a \equiv^{at}_\mcN b$, so part~(ii)
gives $\tp_\mcN(a) = \tp_\mcN(b)$.
If $a, b \in [\mbbN]^2$, then it is clear from the definition of `cl' that $\cl(a)$ and $\cl(b)$ can be ordered as
$a, a', a''$ and $b, b', b''$, respectively, so that $(a, a', a'') \equiv^{at}_\mcN (b, b', b'')$ and again we use part~(ii)
to get $\tp_\mcN(a, a', a'') = \tp_\mcN(b, b', b'')$.
\hfill $\square$

\begin{lem}\label{cl and acl coincide in N}
For every $A \subseteq \mbbN \cup [\mbbN]^2$, $\acl_\mcN(A) = \cl(A) = \dcl_\mcN(A)$.
\end{lem}

\noindent
{\bf Proof.}
By the definition of `$\cl$' it suffices to prove the lemma for finite $A$.
By Lemma~\ref{cl(A) is a subset of dcl(A)}, 
we have $\cl(A) \subseteq \dcl_\mcN(A) \subseteq \acl_\mcN(A)$.
Hence it suffices to show that if $b \notin \cl(A)$ then $b \notin \acl_\mcN(A)$.

Suppose that $b \notin \cl(A)$.
Let $\cl(A) = \{a_1, \ldots, a_n\}$ and let $b'_1, \ldots, b'_m$ enumerate 
$\cl(A \cup \{b\}) \setminus (\cl(A) \cup \{b\})$.
By  Lemma~\ref{extensions of closed sets by elements in F}~(ii) it is enough to find
distinct $b_i$, for $i < \omega$, such that, for each $i < \omega$,
$\cl(A \cup \{b_i\}) \setminus (\cl(A) \cup \{b_i\})$ can be enumerated as
$b'_{i,1}, \ldots, b'_{i,m}$ so that
\[
(a_1, \ldots, a_n, b, b'_1, \ldots, b'_m) \equiv^{at}_\mcN (a_1, \ldots, a_n, b_i, b'_{i,1}, \ldots, b'_{i,m}).
\]
To show this one can argue similarly as in the proof of part~(i) of Lemma~\ref{extensions of closed sets by elements in F}
(hence using Lemma~\ref{N has extension properties}).
The details are left for the reader.
\hfill $\square$

\begin{lem}\label{the subtypes on pairs determine the atomic type of the closure}
Suppose that $a_1, \ldots, a_n, b_1, \ldots, b_n \in [\mbbN]^2$ and
$(a_i, a_j) \equiv_\mcN (b_i, b_j)$ for all $1 \leq i, j \leq n$.
Then $\cl(a_1, \ldots, a_n) \setminus \{a_1, \ldots, a_n\}$ and 
$\cl(b_1, \ldots, b_n) \setminus \{b_1, \ldots, b_n\}$
can be ordered as $a'_1, \ldots, a'_m$ and $b'_1, \ldots, b'_m$, respectively,
so that 
\[
(a_1, \ldots, a_n, a'_1, \ldots, a'_m) \equiv^{at}_\mcN (b_1, \ldots, b_n, b'_1, \ldots, b'_m).
\]
\end{lem}

\noindent
{\bf Proof.}
This is a straightforward consequence of Lemma~\ref{cl(A) is a subset of dcl(A)}.
\hfill $\square$

\begin{lem}\label{the subtypes on pairs determine the type}
Suppose that $a_1, \ldots, a_n, b_1, \ldots, b_n \in [\mbbN]^2$ and
$(a_i, a_j) \equiv_\mcN (b_i, b_j)$ for all $1 \leq i, j \leq n$.
Then $(a_1, \ldots, a_n) \equiv_\mcN (b_1, \ldots, b_n)$.
\end{lem}

\noindent
{\bf Proof.}
Immediate consequence of Lemmas~\ref{extensions of closed sets by elements in F}~(ii)
and~\ref{the subtypes on pairs determine the atomic type of the closure}.
\hfill $\square$

\begin{defin}\label{definition of the relevant substructure of the bipod}{\rm
Let $\mcM$ be a structure with universe $[\mbbN]^2$ and such that, for every
$p = \tp_\mcN(a, b)$ where $a, b \in [\mbbN]^2$ are distinct, $\mcM$ has a relation symbol $R_p$
which is interpreted as the set of realizations of $p$ in $\mcN$. The vocabulary of $\mcM$ has no other 
relation symbols.
}\end{defin}

\begin{lem}\label{the relevant substructure of the bipod is homogeneous}
(i) For all $\bar{a}, \bar{b} \in [\mbbN]^2$ of the same length, 
$\bar{a} \equiv_\mcN \bar{b}$ if and only if $\bar{a} \equiv_\mcM \bar{b}$.\\
(ii) $\mcM$ is homogeneous and has only one complete 1-type over $\es$.\\
(iii) For every $A \subseteq [\mbbN]^2$, $\acl_\mcM(A) = \cl(A) \cap [\mbbN]^2 = \dcl_\mcM(A)$.
\end{lem}

\noindent
{\bf Proof.}
(i) Let $a_1, \ldots, a_n, b_1, \ldots, b_n \in [\mbbN]^2$, $\bar{a} = (a_1, \ldots, a_n)$ and $\bar{b} = (b_1, \ldots, b_n)$.
If $\bar{a} \equiv_\mcM \bar{b}$, then in particular $(a_i, a_j) \equiv_\mcM (b_i, b_j)$ for all $i, j$.
By the definition of $\mcM$ we get $(a_i, a_j) \equiv_\mcN (b_i, b_j)$ for all $i, j$, and then
Lemma~\ref{the subtypes on pairs determine the type}
gives $\bar{a} \equiv_\mcN \bar{b}$.
If $\bar{a} \equiv_\mcN \bar{b}$, then, as $\mcN$ is $\omega$-categorical and countable, there is an
automorphism $\sigma$ of $\mcN$ such that $\sigma(\bar{a}) = \bar{b}$.
Since $[\mbbN]^2$ is $\es$-definable in $\mcN$ (by $\neg F(x)$), $\sigma$ fixes $[\mbbN]^2$ setwise.
From the definition of $\mcM$ it now follows that 
the restriction of $\sigma$ to  $[\mbbN]^2$ is an automorphism of $\mcM$ and hence $\bar{a} \equiv_\mcM \bar{b}$.

Part (ii) follows from~(i) and lemmas~\ref{the subtypes on pairs determine the type}
and~\ref{extensions of closed sets by elements in F}.
Part~(ii) follows from~(i) and
Lemma~\ref{cl and acl coincide in N}
\hfill $\square$

\begin{lem}\label{characterization of dividing in the bipode}
For all tuples $\bar{a}, \bar{b}, \bar{c}$ of elements from $[\mbbN]^2$, the following holds
regardless of whether dividing is considered in $\mcM$ or in $\mcN$:
$\bar{a} \underset{\bar{c}}{\nind} \bar{b}$
if and only if 
\begin{itemize}
\item[(a)] there is $a \in \bar{a}$ such that $a \in \cl(\bar{b}) \setminus \cl(\bar{c})$, or
\item[(b)] there are $a \in \bar{a}$ and $b \in \bar{b}$ such that $a \cap b \neq \es$, but $a \cap c = \es$ for all $c \in \bar{c}$.
\end{itemize}
\end{lem}

\noindent
{\bf Proof sketch.}
If (a) or (b) holds, then it is straightforward to show that $\bar{a} \underset{\bar{c}}{\nind} \bar{b}$ 
(regardless of whether dividing
is considered in $\mcM$ or in $\mcN$).
If neither (a) nor (b) holds, then
one can use
Lemma~\ref{N has extension properties} similarly as in
the proof of 
Lemma~\ref{extensions of closed sets by elements in F}
to show that $\bar{a} \underset{\bar{c}}{\ind} \bar{b}$ (again regardless of whether dividing is in  $\mcM$ or in $\mcN$).
We leave the details to the reader.
\hfill $\square$

\begin{lem}\label{the bipode is simple}
$\mcM$ is supersimple with SU-rank 2, but not stable.
\end{lem}

\noindent
{\bf Proof.}
To prove that $\mcM$ is supersimple it suffices to prove
(by \cite[Theorem~2.4.7 and Definition~2.8.12]{Wag}) 
that if $\mcM' \equiv \mcM$,
$\bar{a} \in M'$ and $B \subseteq M'$, then there is a finite $C \subseteq B$ such that
$\tp_{\mcM'}(\bar{a} / B)$ does not divide over $C$.

Let $\varphi(x, y)$ be a formula in the language of $\mcM$ such that 
for all $a, b \in M (= [\mbbN]^2)$, $\mcM \models \varphi(a, b)$ if and only if $a \neq b$ and $a \cap b \neq \es$.
Recall that $\cl(A) = \dcl_\mcM(A)$ for every $A \subseteq M$.
From Lemma~\ref{characterization of dividing in the bipode} it now follows that
for any $\mcM' \models Th(\mcM)$ and any $\bar{a}, \bar{b}, \bar{c} \in M'$,
$\bar{a} \underset{\bar{c}}{\nind} \bar{b}$ if and only if either
there is some $a \in \bar{a}$ such that $a \in \dcl_{\mcM'}(\bar{b}) \setminus \dcl_{\mcM'}(\bar{c})$,
or there is $b \in \bar{b}$ such that $\mcM' \models \varphi(a, b)$, but $\mcM' \models \neg\varphi(a, c)$ for all $c \in \bar{c}$.

Now suppose that $\mcM' \models Th(\mcM)$, $\bar{a} \in M'$ and $B \subseteq M'$.
Let $C' = \bar{a} \cap \dcl_{\mcM'}(B)$.
For every $a \in \bar{a} \setminus B$ such that there is $b \in B$ such that $\mcM' \models \varphi(a, b)$, 
choose exactly one such $b$ and call it $b_a$.
Let $C = C' \cup \{b_a : a \in \bar{a}\}$. 
Then, for every finite $B' \subseteq B$, $\tp_{\mcM'}(\bar{a} / B'C)$ does not divide over $C$.
By the finite character of dividing, 
$\tp_{\mcM'}(\bar{a} / B)$ does not divide over $C$.

We leave the verification that $\mcM$ has SU rank 2 to the reader.
By using Lemma~\ref{N has extension properties}, it is straightforward to see that
$\mcN$ has the independence property. From this one can derive that also $\mcM$ has the
independence property,  from which it follows that it is unstable.
\hfill $\square$
\\

\noindent
Consider the following equivalence relation on $[\mbbN]^2$:
\[
E_B(a, b) \ \Longleftrightarrow \ \text{ there is $m \in a \cap b$ such that $\mcN \models B(a, m) \wedge B(b, m)$}.
\]
It is clearly $\es$-definable in $\mcN$ and hence it is $\es$-definable in $\mcM$.
By replacing `$B$' with `$R$' we get a similar $\es$-definable equivalence relation $E_R$.
The equivalence classes of $E_B$ and $E_R$ correspond to elements of $\mcM\meq$.
It follows from the definitions of $E_B$, $E_R$ and choice of $\mcN$, that for all $a, b \in [\mbbN]^2$, 
$E_B(a, b) \wedge E_R(a, b)$ if and only if $a = b$.

Let $a \in [\mbbN]^2$. 
By using Lemma~\ref{characterization of dividing in the bipode} and basic ``forking/dividing calculus''
one can now show that, for every $a \in [\mbbN]^2$, $\su(a_{E_B}) = \su(a_{E_{R}}) =~1$ and
$\su(a_= / a_{E_B}) = \su(a_= / a_{E_R}) =~1$ (where clearly $a \in \acl_{\mcM\meq}(a_=)$).
This proves the claim made in
Remark~\ref{remarks on the coordinatization by equivalence relations}~(i), namely that the 
``coordinatization sequence'' of equivalence relations, called $R_1, \ldots, R_k$ in
Theorem~\ref{theorem coordinatizing M by equivalence classes}~(i), need not be unique.

\begin{lem}\label{only two nontrivial equivalence relations in the bipode}
$E_B$ and $E_R$ are the only nontrivial $\es$-definable (in $\mcM$) equivalence relations on $M = [\mbbN]^2$.
\end{lem}

\noindent
{\bf Proof.}
Suppose that $E$ is a nontrivial $\es$-definable (in $\mcM$) equivalence relation on $[\mbbN]^2$
and that $E \neq E_B$ and $E \neq E_R$.
Suppose that $a, b \in [\mbbN]^2$ are such that $a \cap b = \es$ and $E(a, b)$.
Then one can prove,
using Lemma~\ref{N has extension properties}, 
that $E(a', b')$ for all $a', b' \in [\mbbN]^2$, contradicting that $E$ is nontrivial.
We do not give the details, but the idea is that, for any $a', b' \in [\mbbN]^2$,
one can (by Lemma~\ref{N has extension properties}) find $c \in [\mbbN]^2$ such that 
$(a', c) \equiv_\mcN (b', c) \equiv_\mcN (a, b)$, and consequently 
$E(a', c)$ and $E(b', c)$, and thus $E(a', b')$.
Hence, we conclude that, for all $a, b \in [\mbbN]^2$, $E(a, b)$ implies that $a \cap b \neq \es$.

Using the construction of $\mcM$, one can show that there are exactly two binary nontrivial $\es$-definable 
relations which properly refine $E_B$, and none of these two relations is symmetric, hence none of them is 
an equivalence relation. In the same way one can show that there is no nontrivial $\es$-definable equivalence relation
which properly refines $E_R$.
From this (and since $E \neq E_B$ and $E \neq E_R$)  it follows that $E$ does not refine $E_B$ or $E_R$.

Suppose that for all $a, b \in [\mbbN]^2$, $E(a, b)$ implies $E_B(a, b) \vee E_R(a, b)$.
Since $E$ does not refine $E_B$ or $E_R$, and since $\mcM$ has a unique 1-type over $\es$, it follows that there are
distinct $a, b, c \in M$ such that $E(a, b), E(b, c)$, $E_B(a, b)$ and $E_R(b, c)$.
Then $E(a, c)$, so by assumption, $E_B(a, c)$ or $E_R(a, c)$.
But neither case is possible because $E_B(a, b)$ and $E_R(b, c)$.

Hence, there are $a, b \in [\mbbN]^2$ such that 
$E(a, b)$, $\neg E_B(a, b)$ and $\neg E_R(a, b)$ (so $a \neq b$).
Then there is $m \in a \cap b$ such that $\mcN \models B(a, m) \wedge R(b, m)$ or vice versa.
Without loss of generality, suppose that $\mcN \models B(a, m) \wedge R(b, m)$.
Then all $a', b'$ such that $\tp_\mcM(a', b') = \tp_\mcM(a, b)$ or $\tp_\mcM(b', a') = \tp_\mcM(a, b)$
satisfy $E(a', b')$. 
Since 
(by Lemma~\ref{the relevant substructure of the bipod is homogeneous}) 
$\tp_\mcM(a) = \tp_\mcM(b)$, there is
$c \in [\mbbN]^2$ such that $\tp_\mcM(a, b) = \tp_\mcM(b, c)$, so in particular, $E(b, c)$.
Then $\tp_\mcN(a, b) = \tp_\mcN(b, c)$ so there is $n \in b \cap c$ such that
$\mcN \models B(b, n) \wedge R(c, n)$.
Since $\mcN \models R(b, m)$ we have $n \neq m$.
Since $\acl_\mcM(b) = \cl(b) \cap [\mbbN]^2 = \{b\}$ and $b \neq c$ (because $a \neq b$)
we can assume that $c \notin \acl_\mcM(a, b)$, from which it follows (together with $n \in b \cap c$)
that $a \cap c = \es$. But then $\neg E(a, c)$, contradicting the transitivity of $E$.
\hfill $\square$
\\

\noindent
Now we prove the claim made in Remark~\ref{remarks on the coordinatization by equivalence relations}~(ii).
Suppose that $a, b \in [\mbbN]^2$, $a \neq b$, $m \in a \cap b$,
$\mcN \models B(a, m) \wedge R(b, m)$. 
Then 
(by Lemma~\ref{characterization of dividing in the bipode}) 
$a \nind b$, $a_{E_B} \in \acl_\mcM(b)$, 
and
(by some standard forking calculus) 
$a \nind a_{E_B}$.
However, by Lemma~\ref{only two nontrivial equivalence relations in the bipode},
there is {\em no} $\es$-definable equivalence relation $E$ such that $E(a, b)$ and
$a \nind a_E$.

Now we prove Lemma~\ref{the first obstacle} again, this time giving a ``stronger'' example
than in Section~\ref{Cross cutting equivalence relations}
in the sense that, with the notation of Lemma~\ref{the first obstacle}, 
$a \underset{c_R}{\nind} b$ but $b \underset{c_R}{\ind} d$.
By the choice of $\mcN$ and 
Lemma~\ref{N has extension properties},
there are distinct $i, j, k, l, m \in \mbbN$ such that, with
$a = \{i, j\}$, $b = \{k, l\}$, $c = \{j, l\}$ and $d = \{l, m\}$, the following holds in $\mcN$:
\begin{align*}
B(a, j), \ R(c, j), \ B(b, l), \ R(d, l).
\end{align*}
Then $c \underset{c_{R_B}}{\ind} d$, $b \underset{c_{R_B}}{\ind} d$, and $a \underset{c_{R_B}}{\nind} c$.
(The somewhat tedious, but standard, verifications of this are left to the reader.)
Suppose, for a contradiction, that the extension problem (in $\mcM$) of $\tp_\mcM(a, c)$ and $\tp_\mcM(b, d)$ 
has a solution $e = \{i', j'\}$.
Then $i' = j$ or $j' = j$. We can as well assume that $j' = j$.
Since $e \neq c$ we get $i' \neq l$.
As $b \cap d \neq \es$ we must have $e \cap d \neq \es$, which gives $i' = m$.
Hence $e = \{j, m\}$. 
Since $\tp_\mcM(e, c) = \tp_\mcM(a, c)$
we get $\tp_\mcN(e, c) = \tp_\mcN(a, c)$.
Hence $B(e, j)$ and consequently $R(e, m)$.
Then 
\[
\mcN \models \exists x \big( R(e, x) \wedge B(d, x)\big) \ \wedge \ \neg \exists x \big( R(b, x) \wedge B(d, x)\big).
\]
Hence $\tp_\mcN(e, d) \neq \tp_\mcN(b, d)$ and therefore $\tp_\mcM(e, d) \neq \tp_\mcM(b, d)$,
which contradicts that $e$ is a solution to the given extension problem.

\subsection{$\omega$-Pedes}\label{omega-pedes}

In this subsection we outline a proof of Lemma~\ref{the second obstacle}.
The constructions and arguments are similar to, but easier than, those in Section~\ref{Bipedes with bicoloured legs}.
Therefore the proofs of the lemmas that follow are left out.
Let $\mcN = (\mbbN, F^\mcN, E_0^\mcN, E_1^\mcN)$ where:
\begin{itemize}
\item $F$ is unary and $F^\mcN$ and $\mbbN \setminus F^\mcN$ are infinite.
\item $E_0^\mcN$ and $E_1^\mcN$ are equivalence relations such that $E_1^\mcN \subseteq E_0^\mcN$.
\item $E_0$ partitions $F^\mcN$ into infinitely many infinite equivalence classes and $E_1^\mcN$ 
partitions each $E_0^\mcN$-class into exactly two $E_1^\mcN$-classes, both of which are infinite.
\item All $a, b \in \mbbN \setminus F^\mcN$ belong to the same $E_1^\mcN$-class (hence to the same $E_0^\mcN$-class).
\end{itemize}
Let $L$ be a binary relation symbol and let $\Omega$ be the set of expansions 
\[
\mcM = (\mbbN, F^\mcM, E_0^\mcM, E_1^\mcM, L^\mcM)
\]
of $\mcN$ which have the following properties:
\begin{itemize}
\item $\mcM \models \forall x, y \Big( L(x,y) \rightarrow \big(\neg F(x) \wedge F(y)\big)\Big)$.
\item For every $a \in \mbbN \setminus F^\mcN$, every $E_0^\mcN$-class $X \subseteq F^\mcN$ and distinct $E_1^\mcN$-classes
$Y, Z \subseteq X$, either $\mcM \models L(a, b)$ for all $b \in Y$ and $\mcM \models \neg L(a, c)$
for all $c \in Z$, or vice versa.
\end{itemize}
Let 
\[
\Psi = \big\{ (a, X) : a \in \mbbN \setminus F^\mcN \text{ and $X \subseteq F^\mcN$ is an $E_0^\mcN$-class.}\big\}
\]
Let $\Sigma$ and $\mu$ be precisely as in Section~\ref{Bipedes with bicoloured legs}.
Let $\lambda : \Psi \to \mbbN$ be a bijection and 
let $Y_i$, $i < \omega$, be an enumeration of all $E_1^\mcN$-classes which are included in $F^\mcN$.
For every $f \in 2^\mbbN$, let $\mcM_f$ be the unique structure in $\Omega$ which has the following property:
\begin{itemize}
\item[] For every $(a, X) \in \Psi$ and $Y_i, Y_j \subseteq X$, where $i < j$, 
if $f(\lambda(a, X)) = 0$ then $\mcM \models L(a, b)$ for all $b \in Y_i$ and $\mcM \models \neg L(a, c)$ for all $c \in Y_j$,
otherwise $\mcM \models \neg L(a, b)$ for all $b \in Y_i$ and $\mcM \models L(a, c)$ for all $c \in Y_j$.
\end{itemize}

\noindent
Moreover, for every  $\mcM \in \Omega$ there is a unique $f \in 2^\mbbN$ such that $\mcM = \mcM_f$.
In a similar spirit as in the proof of 
Lemma~\ref{extensions of closed sets by elements in F}
(but easier),
one can now prove the following:

\begin{lem}\label{omega-pedes satisfy extension properties}
There is $\mcM \in \Omega$ with the following properties:\\
(i) For all $0 < n < \omega$, all $a_1, \ldots, a_n \in \mbbN \setminus F^\mcM$ and every $f : \{1, \ldots, n\} \to \{0,1\}$,
there is an $E_0^\mcM$-class $X \subseteq F^\mcM$ with $Y_i, Y_j \subseteq X$, where $i < j$, such that
\begin{itemize}
\item[] for every $1 \leq k \leq n$, if $f(k) = 0$ then $\mcM \models L(a_k, b)$ for all $b \in Y_i$ (and hence 
$\mcM \models \neg L(a, c)$ for all $c \in Y_j$), and otherwise $\mcM \models \neg L(a_k, b)$ for all $b \in Y_i$
(and hence $\mcM \models L(a, c)$ for all $c \in Y_j$).
\end{itemize}
(ii) For all $0 < n < \omega$, all $E_0^\mcM$-classes $X_1, \ldots, X_n$ and every $f : \{1, \ldots, n\} \to \{0,1\}$,
there is $a \in \mbbN \setminus F^\mcM$ such that
\begin{itemize}
\item[] for every $1 \leq k \leq n$ and $Y_i, Y_j \subseteq X_k$, where $i < j$,
if $f(k) = 0$ then $\mcM \models L(a, b)$ for every $b \in Y_i$, and otherwise $\mcM \models \neg L(a, b)$ for every $b \in Y_i$.
\end{itemize}
\end{lem}

\noindent
{\bf \em For the rest of this subsection assume that $\mcM$ is like in
Lemma~\ref{omega-pedes satisfy extension properties}.}
Using Lemma~\ref{omega-pedes satisfy extension properties},
one can prove the following by a standard back-and-forth argument which builds up an automorphism:

\begin{lem}\label{the omega-pede is homogeneous}
$\mcM$ is homogeneous.
\end{lem}

\noindent
It is straightforward to see, using 
Lemma~\ref{omega-pedes satisfy extension properties},
that for every $A \subseteq M$, $\acl_\mcM(A) = A$.
With this at hand, it is also straightforward to characterize dividing as follows:

\begin{lem}\label{dividing in the omega-pede}
For all $\bar{a}, \bar{b}, \bar{c} \in M$, $\tp_\mcM(\bar{a} / \bar{b}\bar{c})$ divides over $\bar{c}$
if and only if there is $a \in \bar{a}$ such that
\begin{itemize}
\item[(i)] $a \in \bar{b} \setminus \bar{c}$, or
\item[(ii)] $\mcM \models F(a)$ and there is $b \in \bar{b}$ such that $\mcM \models E_0(a, b)$ and for all
$c \in \bar{c}$, $\mcM \models \neg E_0(a, c)$.
\end{itemize}
\end{lem}

\noindent
With Lemma~\ref{dividing in the omega-pede} and standard
arguments as in the proof of Lemma~\ref{the bipode is simple},
one can prove:

\begin{lem}\label{the omega-pede is simple}
$\mcM$ is supersimple (but not stable). If $\mcM \models F(a)$ then $\su(a) = 2$, otherwise $\su(a) = 1$.
\end{lem}

\noindent
Now we are ready to prove Lemma~\ref{the second obstacle}.
There are $c, d \in M$ such that
\[
\mcM \models F(c) \wedge F(d) \wedge E_0(c, d) \wedge \neg E_1(c, d).
\]
By Lemma~\ref{omega-pedes satisfy extension properties}, we can also find $a, b \in M$ such that 
$\mcM \models L(a, c) \wedge L(b, d)$, and hence $\mcM \models \neg F(a) \wedge \neg F(b)$.
Since $\mcM$ is homogeneous there is an automorphism of $\mcM\meq$ which
takes $(a, c)$ to $(b, d)$. 
This automorphism can be extended to an automorphism of $\mcM\meq$.
Since $E_0(c, d)$ it follows that this automorphism (of $\mcM\meq$) fixes $c_{E_0}$.
Hence $\tp_{\mcM\meq}(a / c_{E_0}) = \tp_{\mcM\meq}(b / c_{E_0})$.
But there is no $e$ such that $\tp_\mcM(e, c) = \tp_\mcM(a, c)$ and $\tp_\mcM(e, d) = \tp_\mcM(b, d)$,
because this would give $L(e, c) \wedge L(e, d)$ where $E_0(c, d) \wedge \neg E_1(c, d)$.
However note that $\tp_{\mcM\meq}(a / \acl_{\mcM\meq}(c_{E_0})) \neq \tp_{\mcM\meq}(b / \acl_{\mcM\meq}(c_{E_0}))$,
because $c_{E_1}, d_{E_1} \in \acl_{\mcM\meq}(c_{E_0})$.

\subsection{Metric spaces}\label{Metric spaces}\footnote{
I thank G. Conant for explaining how results in \cite{Con-neostab} relate to results in this article.}
Unlike sections~\ref{Cross cutting equivalence relations} --~\ref{omega-pedes},
the examples of this section
are not meant to show that things can be more complicated than one might have hoped for.
Instead these are examples for which the main results of this article
are concretized in nice and natural ways.

In \cite{Con-neostab}, Conant studies the infinite countable homogeneous (in a more general sense that in this article)
metric space, denoted $\mcU_\mcR$ and called {\em $\mcR$-Urysohn space}, 
over a countable {\em distance monoid} $\mcR = (R, \oplus, \leq, 0)$
(see \cite[Section~2]{Con-neostab} for a definition of distance monoid).
In other words,
fix some distance monoid $\mcR$ and let $\mcK_\mcR$
be the class of all finite $\mcR$-metric spaces.
Then, for a suitable relational language, $\mcK_\mcR$ is closed under isomorphism
and has the hereditary property and the amalgamation property.
Hence the Fra\"{i}ss\'{e} limit of $\mcK_\mcR$ exists and we denote it by $\mcU_\mcR$.
The language that we use has a binary relation symbol
$d_r$ for every $r \in R$, where $d_r(a, b)$ is interpreted as ``the distance between $a$ and $b$ is at most $r$''.
So a structure $\mcM$ for this vocabulary is viewed as an {\em $\mcR$-metric space} if for all $a, b, c \in M$,
\begin{itemize}
\item $d_0(a, b)$ if and only if $a = b$,
\item for all $r \in R$, $d_r(a, b)$ if and only if $d_r(b, a)$, and 
\item (triangle inequality) for all $r, s, t \in R$, if $d_r(a, b)$, $d_s(b, c)$ and $d_t(a, c)$, then \\ $r \oplus s \geq t$.
\end{itemize}

\noindent
If $\mcR$ is finite then the vocabulary of an $\mcR$-metric space is finite and hence $\mcU_\mcR$
is homogeneous in the sense of this article.
From now on, assume that $\mcR$ is a {\em finite} distance monoid.
As examples of finite distance monoids one can take $\mcR = (R, \oplus, \leq, 0)$,
where $R \subseteq \mbbR^{\geq 0}$ is finite and chosen so that $0 \in R$, `$\leq$' is the
usual order on $\mbbR$, `$\oplus$' is `$+_R$' where for all $r, s \in R$,
\[
r +_R s = \max\{ x \in R : x \leq r + s\} \quad  \text{ and \ \ $+_R$ is associative.}
\]
For example, this holds if one takes $R = \{0, 1, 2\}$ or $R = \{0, 1, 3, 4\}$.
In the first case, however, $\mcU_\mcR$ is essentially the Rado graph, 
by viewing ``$d_1(a, b) \wedge \neg d_0(a, b)$'' as saying that there is an edge between $a$ and $b$,
and ``$d_2(a, b) \wedge \neg d_1(a, b)$'' as saying that there is no edge between $a$ and $b$ (and $a \neq b$).
More examples of finite distance sets are analyzed in Appendix~A of L. Nguyen van Th\'{e}'s thesis~\cite{vanThe}.

By \cite[Theorem~4.9]{Con-neostab}, {\em $\mcU_\mcR$ is simple if and only if for all $r, s \in R$ such that $r \leq s$,
$r \oplus r \oplus s = r \oplus s$.}
One can check that if, for example, $R = \{0, 1, 3, 4\}$ then this condition holds.
From now on, suppose that $\mcU_\mcR$ is simple. 
Hence it is (by Fact~\ref{facts about binary simple homogeneous structures})
supersimple with finite SU-rank and trivial dependence.
An element $r \in R$ is called {\em idempotent} if $r \oplus r = r$.
By \cite[Theorem~4.16]{Con-neostab}, 
{\em the SU-rank of $\mcU_\mcR$ is the number of non-maximal idempotent elements in $R$.}
Moreover, by \cite[Corollary~7.9]{Con-neostab},
{\em the $\es$-definable equivalence relations on the universe of $\mcU_\mcR$ are exactly those which are defined by
the formulas $d_r(x, y)$ where $r$ is idempotent.}
Suppose that $0 < r < s \in R$ are idempotent elements.
Using the idempotency one can easily show that the equivalence relation $d_r(x, y)$ partitions every class of the 
equivalence relation $d_s(x, y)$ into infinitely many parts, all of which are infinite.
Thus the sequence of equivalence relations $R_1, \ldots, R_k$ in 
Theorem~\ref{theorem coordinatizing M by equivalence classes}~(i) corresponds, in the case of $\mcU_\mcR$,
to $d_{r_1}(x, y), \ldots, d_{r_k}(x, y)$, where $r_1 > \ldots > r_k$ is a list of all non-maximal idempotent elements
(so $r_k = 0$).

For any $r \in R$, let `$2r$' denote `$r \oplus r$'.
From the characterization of $\mcU_\mcR$ being simple (given above), it follows that $2r$ is idempotent for every $r \in R$.
Let $d(a, b)$ be the least $r \in R$ such that $d_r(a, b)$ holds.
From \cite[Corollary~4.10]{Con-neostab} we have for all $a, b$ and $\bar{c}$
from any model of $Th(\mcU_\mcR)$:
\[
a \underset{\bar{c}}{\nind} b \ \Longleftrightarrow \ 2d(a, b) \ < \ 2d(a, c) \ \text{ for all $c \in \bar{c}$.}
\]
Since $2r$ is idempotent for every $r \in R$, it follows that, for every $r \in R$,
$a \underset{\bar{c}}{\nind} b$ if and only if there is a $\es$-definable equivalence relation $E$,
defined by $d_r(x, y)$ for some idempotent $r$, such that $E(a, b)$ but $\neg E(a, c)$ for all $c \in \bar{c}$.
This is the specific version of Theorem~\ref{theorem coordinatizing M by equivalence classes}~(iii)
in the case of $\mcU_\mcR$.


\begin{thebibliography}{99}\label{References}

\bibitem{Ahl} O. Ahlman, Simple theories axiomatized by almost sure theories,
{\em Annals of Pure and Applied Logic}, Vol. 167 (2016) 435--456.

\bibitem{AK} O. Ahlman, V. Koponen, On sets with rank one in simple homogeneous structures,
{\em Fundamenta Mathematicae}, Vol. 228 (2015) 223--250.

\bibitem{AL} A. Aranda L\'{o}pez, {\em Omega-categorical simple theories}, Ph.D. thesis, 
The University of Leeds (2014).

\bibitem{Che98} G. L. Cherlin, {\em The Classification of Countable Homogeneous Directed Graphs
and Countable Homogeneous $n$-tournaments}, 
Memoirs of the American Mathematical Society 621, American Mathematical Society (1998).

\bibitem{CL} G. Cherlin, A. H. Lachlan, Stable finitely homogeneous structures,
{\em Transactions of the American Mathematical Society}, Vol. 296 (1986) 815--850.

\bibitem{Con} G. Conant, An axiomatic approach to free amalgamation, {\em The Journal of Symbolic Logic}, to appear. Available at: 
\url{http://arxiv.org/abs/1505.00762}

\bibitem{Con-neostab} G. Conant, Neostability in countable homogeneous metric spaces.
Available at: 
\url{https://arxiv.org/abs/1504.02427}

\bibitem{DK} T. De Piro, B. Kim, The geometry of 1-based minimal types,
{\em Transactions of The American Mathematical Society}, Vol. 355 (2003) 4241--4263.

\bibitem{Djo06} M. Djordjevi\'{c}, Finite satisfiability and $\omega$-categorical structures with trivial dependence,
{\em The Journal of Symbolic Logic}, Vol. 71 (2006) 810--829.

\bibitem{Fra54} R. Fra\"{i}ss\'{e}, Sur l'extension aux relations de quelques propri\'{e}t\'{e}s des ordres,
{\em Annales Scientifiques de l'\'{E}cole Normale Sup\'{e}rieure}, Vol. 71 (1954) 363--388.

\bibitem{Fri} A. Friedman, {\em Foundations of Modern Analysis}, Dover Publications, New York (1982).

\bibitem{Gar} A. Gardiner, Homogeneous graphs, 
{\em Journal of Combinatorial Theory, Series B}, Vol. 20 (1976) 94--102.

\bibitem{GK} Y. Golfand, M. Klin, On $k$-homogeneous graphs, 
in {\em Algorithmic Studies in Combinatorics}, Nauka, Moscow (1978), 76--85.

\bibitem{Goode} J. B. Goode, Some trivial considerations, {\em The Journal of Symbolic Logic},
Vol. 56 (1991) 624--631.

\bibitem{HKP} B. Hart, B. Kim, A. Pillay, Coordinatisation and canonical bases in simple theories,
{\em The Journal of Symbolic Logic}, Vol. 65 (2000) 293--309.

\bibitem{Hod} W. Hodges,  {\em Model theory}, Cambridge University Press (1993).

\bibitem{KL} J. Knight, A. H. Lachlan, Shrinking, stretching and codes for homogeneous structures,
{\em Classification Theory}, Lecture Notes in Mathematics 1292, Springer Verlag, Berlin--New York, 192--228 (1987).

\bibitem{Kop09} V. Koponen, Independence and the finite submodel property,
{\em Annals of Pure and Applied Logic}, Vol. 158 (2009) 58--79.

\bibitem{Kop16a} V. Koponen, Binary simple homogeneous structures are supersimple with finite rank,
{\em Proceedings of the American Mathematical Society}, Vol. 144 (2016) 1745--1759.

\bibitem{Kop-one-based} V. Koponen, Homogeneous 1-based structures and interpretability in random structures, 
{\em Mathematical Logic Quarterly}, to appear.

\bibitem{Kop16b} V. Koponen, Binary primitive homogeneous simple structures, 
{\em The Journal of Symbolic Logic}, to appear.

\bibitem{Lach84} A. H. Lachlan, Countable homogeneous tournaments,
{\em Transactions of the American Mathematical Society}, Vol. 284 (1984) 431--461.

\bibitem{Lach97} A. H. Lachlan, Stable finitely homogeneous structures: a survey,
in B. T. Hart et. al. (eds.), {\em Algebraic Model Theory}, 145--159, Kluwer Academic Publishers (1997) 

\bibitem{LS} A. H. Lachlan, S. Shelah, Stable structures homogeneous for a finite binary language,
{\em Israel Journal of Mathematics}, Vol. 49 (1984) 155--180.

\bibitem{LachTripp} A. H. Lachlan, A. Tripp, Finite homogeneous 3-graphs, 
{\em Mathematical Logic Quarterly}, Vol. 41 (1995) 287--306.

\bibitem{LW} A. H. Lachlan, R. Woodrow,  Countable ultrahomogenous undirected graphs,
{\em Transactions of the Americal Mathematical Society}, Vol. 262 (1980) 51--94.

\bibitem{LT} D. Lockett, J. K. Truss, Homogeneous coloured multipartite graphs,
{\em European Journal of Combinatorics}, Vol. 42 (2014) 217--242.

\bibitem{Mac91} D. Macpherson, Interpreting groups in $\omega$-categorical structures,
{\em The Journal of Symbolic Logic}, Vol. 56 (1991) 1317--1324.

\bibitem{Mac11} D. Macpherson, A survey of homogeneous structures,
{\em Discrete Mathematics}, Vol. 311 (2011) 1599--1634.

\bibitem{Schm} J. H. Schmerl, Countable homogeneous partially ordered sets,
{\em Algebra Universalis}, Vol. 9 (1979) 317--321.

\bibitem{Shee} J. Sheehan, Smoothly embeddable subgraphs,
{\em Journal of The London Mathematical Society}, Vol. 9 (1974) 212--218.

\bibitem{She} S. Shelah, {\em Classification Theory}, Revised Edition, North-Holland (1990).

\bibitem{vanThe} L. Nguyen van Th\'{e},
{\em Structural Ramsey Theory of Metric Spaces and Topological Dynamics of Isometry Groups},
Memoirs of the American Mathematical Society 968, Americal Mathematical Society, Providence (2010).

\bibitem{Wag} F. O. Wagner, {\em Simple Theories}, Kluwer Academic Publishers (2000).

\end{thebibliography}
\end{document}